\documentclass[leqno,11pt]{amsart}

\usepackage{geometry}

\usepackage[all]{xy} 

\usepackage{comment}

\usepackage{mathtools}

\usepackage{amsmath, amssymb, amsfonts, latexsym, mdwlist, amsthm}
\usepackage{subfig}
\usepackage{graphicx}
\usepackage{wrapfig}

\usepackage[bookmarks, colorlinks, breaklinks, pdftitle={},
pdfauthor={Soheyla Feyzbakhsh}]{hyperref}
\hypersetup{linkcolor=blue,citecolor=blue,filecolor=black,urlcolor=blue}

\usepackage{tikz}
\usetikzlibrary{calc,trees,positioning,arrows,chains,shapes.geometric,%
    decorations.pathreplacing,decorations.pathmorphing,shapes,%
    matrix,shapes.symbols}

\tikzset{
>=stealth',
  punktchain/.style={
    rectangle,
    rounded corners,
    draw=black, thick,
    minimum height=3em,
    text centered,
    on chain},
  line/.style={draw, thick, <-},
  element/.style={
    tape,
    top color=white,
    bottom color=blue!50!black!60!,
    minimum width=8em,
    draw=blue!40!black!90, very thick,
    text width=10em,
    minimum height=3.5em,
    text centered,
    on chain},
  every join/.style={->, thick,shorten >=1pt},
  decoration={brace},
  tuborg/.style={decorate},
  tubnode/.style={midway, right=2pt},
}

\usepackage{paralist}
\setdefaultenum{(a)}{(i)}{}{}
\usepackage{enumitem} 

\renewcommand\_{^{}_}

\newcommand\ch{\operatorname{ch}}
\newcommand\coh{\operatorname{Coh}}
\newcommand\rk{\operatorname{rk}}
\newcommand\cO{\mathcal O}

\newcommand\al{\alpha}

\newcommand\cH{\mathcal H}

\renewcommand\_{^{}_}
\renewcommand\;{\hspace{.6pt}}

\newcommand\PP{\mathbb P}

\newcommand\Z{\mathbb Z}

\newcommand\cA{\mathcal A}
\newcommand\cD{\mathcal D}
\newcommand\Pic{\operatorname{Pic}}
\newcommand\Cliff{\operatorname{Cliff}}

\newcommand\Coh{\operatorname{Coh}}
\newcommand\cok{\operatorname{coker}}

\newcommand\im{\operatorname{im}}

\def\abs#1{\left\lvert#1\right\rvert}

\makeatletter
\newtheorem*{rep@theorem}{\rep@title}
\newcommand{\newreptheorem}[2]{%
\newenvironment{rep#1}[1]{%
 \def\rep@title{#2 \ref{##1}}%
 \begin{rep@theorem}}%
 {\end{rep@theorem}}}
\makeatother

\newtheorem{Thm}{Theorem}[section]
\newreptheorem{Thm}{Theorem}
\newtheorem{Prop}[Thm]{Proposition}

\newtheorem{Lem}[Thm]{Lemma}

\newtheorem{Cor}[Thm]{Corollary}
\newreptheorem{Cor}{Corollary}

\newreptheorem{Con}{Conjecture}

\newtheorem{thm-int}{Theorem}

\theoremstyle{definition}
\newtheorem{Def-s}[Thm]{Definition}
\newtheorem{Def}[Thm]{Definition}
\newtheorem{Rem}[Thm]{Remark}

\newcommand{\ignore}[1]{}

\begin{document}


\title[An effective restriction theorem]
{An effective restriction theorem via Wall-crossing\\ and Mercat's conjecture}

\author{Soheyla Feyzbakhsh}
\address{Department of Mathematics,
Imperial College,
London SW7 2AZ, United Kingdom}
\email{s.feyzbakhsh@imperial.ac.uk}


\begin{abstract}
We prove an effective restriction theorem for stable vector bundles $E$ on a smooth projective variety: $E|_D$ is (semi)stable for all irreducible divisors $D \in |kH|$ for all $k$ greater than an explicit constant. 

As an application, we show that Mercat's conjecture in any rank greater than $2$ 
fails for curves lying on K3 surfaces. 

Our technique is to use wall-crossing with respect to (weak) Bridgeland stability conditions which we also use to reprove Camere's result on slope stability of the tangent bundle of $\mathbb{P}^n$ restricted to a K3 surface.  
 
\end{abstract}

\vspace{-3cm}
\maketitle

\vspace{-.4cm}
\section{Introduction}

Inspired by the construction of Bridgeland stability conditions on K3 surfaces \cite{bridgeland:K3-surfaces} Toda, Bayer, Macr\`i and Stellari \cite{bayer:bridgeland-stability-conditions-on-threefolds,bayer:the-space-of-stability-conditions-on-abelian-threefolds} studied weak stability conditions on any smooth complex projective variety. In this paper, we use wall-crossing with respect to these weak stability conditions to prove an effective restriction theorem that expresses sufficient conditions on a slope-stable reflexive sheaf such that its restriction to a hypersurface remains stable. Restriction theorems provide us with the possibility of studying higher dimensional varieties via the geometry of their subvarieties. That is why they have been long-studied via different approaches; see \cite[Chapter 7]{huybrechts:geometry-of-moduli-space-of-sheaves} for a survey.

Let $X$ be a smooth complex projective variety of dimension $n \geq 2$ with an ample divisor $H$. For a $\mu$-stable coherent sheaf $E$ of positive rank on $X$, we define $$
\widetilde{\Delta}(E) \coloneqq \left(\frac{\ch_1(E).H^{n-1}}{\ch_0(E)H^n}\right)^2 - 2 \frac{\ch_2(E).H^{n-2}}{\ch_0(E)H^n},$$
\begin{align}
&\mu^{\max}(E) \coloneqq \max\left\{\mu(F) \colon \text{$F$ is a subsheaf of $E$ with $\mu(F) <\mu(E)$}\right\}, \label{max} \\
&\mu^{\min}(E) \coloneqq \min\left\{\mu(F') \colon \text{$F'$ is a proper quotient sheaf of $E$ 
}\right\}, \label{min}
\end{align}
and $\delta(E) \coloneqq \min\left\{\mu^{\min}(E) - \mu(E),\, \mu(E) -\mu^{\max}(E)  \right\}$.  
\begin{Thm}\label{thm.main thm}
	Let $E$ be a $\mu$-stable reflexive sheaf on $X$ of rank $\rk > 0$. The restricted sheaf $E|_D$ for any irreducible divisor\footnote{$D$ can be singular.} $D \in |k H|$ is $\mu$-semistable on $D$ if
	\begin{equation}\label{bound for ell-new}
    k \,\geq\, \frac{\rk +2}{\sqrt{\rk+1}}\sqrt{\widetilde{\Delta}(E)} \qquad \text{and} \qquad   
    \frac{k}{2} + \sqrt{\frac{k^2}{4} - \widetilde{\Delta}(E)} \,\geq\, \frac{\widetilde{\Delta}(E)}{\delta(E)}\ .	 
	\end{equation}
	Moreover, $E|_D$ is $\mu$-stable on $D$ if the inequalities in \eqref{bound for ell-new} are both strict. 
\end{Thm}
The $\mu$-slope of a coherent sheaf and the notion of $\mu$-(semi)stability are defined in section \ref{section.2.background}. Note that if $k \geq 2\delta(E)$, the conditions in \eqref{bound for ell-new} are equivalent to
\begin{equation}\label{bound for ell}
	k \geq \max\left\{ \frac{\rk +2}{\sqrt{\rk+1}}\sqrt{\widetilde{\Delta}(E)}\ , \ \delta(E) + \frac{\widetilde{\Delta}(E)}{\delta(E)}  \right\}. 
	\end{equation}
When $\rk> 1$, we always have $\delta(E) \geq \frac{1}{H^n\rk(\rk -1)} $.
If we substitute this lower bound, 
we obtain one of Langer's restriction theorems \cite
{langer:semistable-sheaves-in-positive-characteristics}, see Corollary \ref{cor-langer} and Remark \ref{rem-langer} for more details.  

\subsection*{Clifford indices} 
The Clifford index Cliff$(C)$ of a smooth curve $C$ is the second important invariant of $C$ after the genus $g$, which carries the information of special line bundles $C$. 
Lange and Newstead \cite{lange:cliiford-indices-for-vector-bundles} proposed a generalisation of Cliff$(C)$ to higher rank Clifford index Cliff$_r(C)$ which depends on rank $r$-semistable vector bundles on $C$. 

Take a vector bundle $E$ of rank $r$ and degree $d$ on $C$. The Clifford index of $E$ is defined as Cliff$(E) = \frac{1}{r}(d-2(h^0(E)-r))$.  We say $E$ contributes to the rank $r$-Clifford index of $C$ if $E$ is $\mu$-semistable with degree $d \leq r(g-1)$ and $h^0(C,E) \geq 2r$. Then the rank $r$-Clifford index of $C$ is defined as the quantity
\begin{equation*}
\Cliff_r(C) \coloneqq \min \big\{\!\Cliff(E) \colon \text{$E$ contributes to the rank $r$-Clifford index of $C$}    \big\}.
\end{equation*}
Note that Cliff$_1(C) = \text{Cliff}(C)$ is the classical Clifford index of $C$. In terms of these new invariants, Mercat's conjecture \cite{mercat:clifford-theorem} can be expressed as 
$$
\text{${\bf M_r}(C) \colon $}  
\;\;\text{Cliff}_r(C) = \text{Cliff}(C)
$$
which says higher rank Clifford indices are equal to the rank one Clifford index.

Curves over K3 surfaces have played an important role in the Brill-Noether theory of vector bundles on curves. In  \cite{bakker:mercat-conjecture}, the conjecture ${\bf M_2}(C)$ was proved for any smooth curve $C \in |H|$ on a K3 surface $S$ when $\text{Pic}(S) = \mathbb{Z}.H$. Using this, ${\bf M_2}(C)$ was shown for generic curves of every genus. However, the restriction of Lazarsfeld-Mukai bundles on $S$ to the curve $C$ (see Section \ref{section.clifford} for a definition) have led to counterexamples to ${\bf M_3}(C)$ \cite{farkas:higher-rank-brill-noether} and ${\bf M_4}(C)$ \cite{aprodu:restricted-lazarsfeld-mukai}. As a consequence of Theorem \ref{thm.k3 surface.cor}, we generalise these results to higher ranks and show ${\bf M_r}(C)$ fails for $r \geq 3$ and any smooth curve $C \in |H|$.
\begin{Thm}\label{thm.k3 surface.cor}
		Let $(S, H)$ be a smooth polarised K3 surface such that 
		\begin{equation}\label{assumption}
		\text{$H^2$ divides $H.D$ for all curve classes $D$ on $S$}, 
		\end{equation}
		e.g. $\Pic(S) = \mathbb{Z}.H$. Take a $\mu$-stable vector bundle $E$ on $S$ with Chern character $\ch(E) = (r, H, \ch_2)$ such that $r \geq 2$ and $\ch_2 \geq 0$. Then for any smooth curve $C \in |H|$ of genus $g$, the restricted bundle $E|_C$ contributes to the rank $r$-Clifford index of $C$. If $r \geq 4$, or, $r=3$ and $g -\ch_2 < 4+ \frac{3}{2}\lfloor \frac{g-1}{2} \rfloor$, then
		\begin{equation*}
		\Cliff(E|_C) <  \Cliff(C) = \left\lfloor \frac{g-1}{2} \right\rfloor .  
		\end{equation*} 
		In particular, ${\bf M_r}(C)$ does not hold if either (i) $r \geq 4$ and $g \geq r^2$, or (ii) $r=3$ and $g=7, 9$ or $g \geq 11$.  
		
\end{Thm}     
\subsection*{Slope stability of the restriction of tangent bundle of $\mathbb{P}^n$ to a K3 surface.}
Let $S$ be a smooth projective K3 surface and let $L$ be an ample line bundle generated by its global sections on $S$. Then there is a well-defined morphism 
\begin{equation*}
\phi_L \colon S \rightarrow \mathbb{P}(H^0(L)^*) \cong \mathbb{P}^n. 
\end{equation*}
In the final part of the paper, we apply wall-crossing with respect to Bridgeland stability conditions on K3 surfaces to reprove Camere's result on slope-stability of $\phi_L^{*}T_{\mathbb{P}^n}$. The restriction of the Euler exact sequence to the K3 surface $S$ and tensoring with $L^*$ gives the short exact sequence  
\begin{equation*}
0 \rightarrow \big(\phi_L^*T_{\PP^n}\big)^* \otimes L \rightarrow H^0(S,L) \otimes \mathcal{O}_S \xrightarrow{\text{ev}} L \rightarrow 0\ . 
\end{equation*}


\begin{Thm}\cite[Theorem 1]{camere:tangent-bundle-of-pn} \label{thm. Theorem of ML}
	Let $S$ be a smooth projective K3 surface over $\mathbb{C}$, and let $L$ be a globally generated ample line bundle on $S$. Then the kernel $M_L$ of the evaluation map on the global sections of $L$
	\begin{equation}\label{def M_L}
	0 \rightarrow M_L \rightarrow H^0(S,L) \otimes \mathcal{O}_S \xrightarrow{\text{ev}} L \rightarrow 0, 
	\end{equation}   
	is $\mu$-stable with respect to $L$.
\end{Thm}

\subsection*{Method of the proof of Theorem \ref{thm.main thm}
}
There is an abelian subcategory $\mathcal{A} \subset \cD(X)$ of the bounded derived category of coherent sheaves on $X$ which includes $E$ and $E(-k H)[1]$ for $k \in \mathbb{N}$. Thus for any irreducible divisor $D \in |k H|$, the restricted sheaf $E|_D$ lies in an exact sequence
\begin{equation*}
E \rightarrow E|_D \rightarrow E(-k H)[1]
\end{equation*}
in $\mathcal{A}$. The slope-stability of the reflexive sheaf $E$ implies that there is a weak stability condition $\sigma$ on $\cA$ such that both $E$ and $E(-k H)[1]$ are stable with respect to $\sigma$. If $k$ is large enough, we may apply wall-crossing techniques to show that we can deform the weak stability condition $\sigma$, while keeping $E$ and $E(-k H)[1]$ stable, to reach a weak stability condition $\sigma'$ where $E$ and $E(-k H)[1]$ have the same slope (phase). Thus their extension $E|_D$ is $\sigma'$-semistable of the same slope. Then a general argument immediately implies that $E|_D$ is slope-stable. 

The main advantage of this method is the possibility to strengthen effective restriction theorems for special sheaves $E$ as soon as we have more control over the position of the walls for $E$ and $E(-kH)[1]$; see for instance Proposition \ref{special}. 


        
\subsection*{Related work}

Theorem \ref{thm.main thm} appeared in an earlier version of this paper \cite{feyz:slope-stability-of-restrcition}, where it was stated for K3 surfaces only. However the proof works for any variety, as the current version of the paper makes explicit. The result has now found applications in \cite{chunyi:stability-condition-quintic-threefold,koseki:double-triple-solids} to construct Bridgeland stability conditions on Calabi-Yau threefolds, and in \cite{kopper:stability-conditions-for-restriction} to investigate the stability of the restriction of stable vector bundles on a surface to any integral curve, not necessarily multiples of $H$.

%
%
%
%
%
%
%
%
%
%
%

In joint work with Chunyi Li \cite{feyz-li:clifford-indices}, we employed the techniques in the current paper and \cite{feyz:mukai-program} to compute rank $\geq 2$ Clifford indices of curves over K3 surfaces. 

\subsection*{Acknowledgements}I would like to thank Arend Bayer for many useful discussions. I am grateful for comments by Gavril Farkas, Chunyi Li, Angela Ortega and Richard Thomas. The author was supported by the ERC starting grant (PI: Arend Bayer) WallXBirGeom 337039 and EPSRC postdoctoral fellowship EP/T018658/1. 

\section{Two-dimensional slice of weak stability conditions}\label{section.2.background}
In this section, we recall the notion of (weak) Bridgeland stability conditions on the bounded derived category of coherent sheaves. The main references are \cite{bayer:the-space-of-stability-conditions-on-abelian-threefolds,bayer:bridgeland-stability-conditions-on-threefolds}. \par
Let $X$ be a projective scheme over $\mathbb{C}$ and let $H$ be the class of an ample divisor on $X$. Recall that the Hilbert polynomial  $P(E, m)$ of a coherent sheaf $E$ is defined via the Euler characteristic $m \mapsto \chi(E \otimes \mathcal{O}_X(mH))$. It can be uniquely written as   
\begin{equation*}
P(E, m) = \sum_{i=0}^{\dim E} \alpha_i(E) \ \frac{m^i}{i !}
\end{equation*}
with integral coefficients $\alpha_i(E)$ for $i=0, \dots , \dim E$ where $\dim E$ is the dimension of the support of $E$ \cite{huybrechts:geometry-of-moduli-space-of-sheaves}. If $E$ is a coherent sheaf of dimension $d = \dim X$, we define the rank and degree of $E$ as 
\begin{equation*}
\rk(E) \coloneqq \frac{\alpha_d(E)}{\alpha_d(\cO_X)} \  ,  \qquad \deg (E) \coloneqq \alpha_{d-1}(E) - \rk(E) \cdot \alpha_{d-1}(\cO_X) .
\end{equation*} 
Then for any coherent sheaf $E$, the slope is defined as  
\begin{equation}\label{slope function}
\mu(E) \coloneqq \ \left\{\!\!\begin{array}{cc} \frac{\deg(E)}{\rk(E)} & \text{if } \alpha_{\dim(X)} \neq 0 
, \\
+\infty & \text{oth.}
\end{array}\right.
\end{equation}
\begin{Def}\label{def:slope-stability}
	A coherent sheaf $E$ of dimension $d = \dim (X)$ is $\mu$-(semi)stable if for all 
    non-trivial quotient sheaves $E \twoheadrightarrow F$ we have $\mu(E) < (\leq) \ \mu(F)$.  
\end{Def}

From now on, we assume $X$ is a smooth projective complex variety of dimension $n \geq 2$. For any coherent sheaf $E$ on $X$, the Hirzebruch-Riemann-Roch formula shows that $\deg(E) = \ch_1(E).H^{n-1}$. 

We denote the bounded derived category of coherent sheaves on $X$ by $\mathcal{D}(X)$ and its Grothendieck group by $K(X)$. 
For any $b \in \mathbb{R}$, consider the full abelian subcategory of complexes 
\begin{equation}\label{Abdef}
\Coh^b(X)\ \coloneqq \ \big\{E^{-1} \xrightarrow{\,d\,} E^0 \ \colon\ \mu^{+}(\ker d) \leq b \,,\  \mu^{-}(\cok d) > b \big\}. 
\end{equation}
Here $\mu_H^{\pm}(E)$ for a coherent sheaf $E$ is the maximum (minimum) slope in the Harder-Narasimhan filtration of $E$ with respect to $\mu$-stability. By \cite[Lemma 6.1]{bridgeland:K3-surfaces}, $\Coh^b(X)$ is the heart of a bounded t-structure on $\cD(X)$. 

For any $w > \frac{b^2}{2}$, we define the function $Z_{b,w} \colon K(X) \rightarrow \mathbb{C}$ such that for any $[E] \in K(X)$, 
\begin{align*}
Z_{b,w}(E) = - \ch_2(E).H^{n-2} + w\ch_0(E)H^n + i\big(\ch_1(E).H^{n-1}-bH^n\ch_0(E)\big).  
\end{align*}
Each $E\in \Coh^b(X)$ satisfies Im$[Z_{b,w}(E)] \geq 0$, and if Im$[Z_{b,w}(E)] = 0$, then Re$[Z_{b,w}(E)] \leq 0$ \cite[Proposition 12.2]{bayer:the-space-of-stability-conditions-on-abelian-threefolds} \footnote{To have a weak stability function whose imaginary part is non-negative and the real part is linear, we replace $\overline{Z}_{a, B = b H}$ in \cite[Proposition 12.2]{bayer:the-space-of-stability-conditions-on-abelian-threefolds} by 
\vspace{-.1 cm}
{\small
\begin{equation*}
i\, \overline{Z}_{\alpha, b H} - \beta\, \text{Im}[i\, \overline{Z}_{\al, b H}] = -\ch_2.H^{n-2} + \ch_0H^n \left(\frac{\al^2}{2} + \frac{b^2}{2} \right)+ i(\ch_1.H^{n-1} -b\ch_0H^n) 
\end{equation*}	
}
\vspace{-.1 cm}	
and substitute $\frac{\alpha^2}{2} + \frac{b^2}{2}$ with $w$. 
}. Thus for objects in $\coh^b(X)$, we have the slope function 
%
%
%
\begin{equation}\label{noo}
\nu\_{b,w}(E)\ \coloneqq\ \left\{\!\!\!\begin{array}{cc} -\frac{\text{Re}[Z_{b,w}(E)]}{\text{Im}[Z_{b,w}(E)]}
& \text{if }\text{Im}[Z_{b,w}(E)]\ne0, \\
+\infty & \text{if }\text{Im}[Z_{b,w}(E)]=0. \end{array}\right.
\end{equation}
\begin{Def}
	Fix $w>\frac{b^2}2$. We say $E\in\cD(X)$ is $\nu\_{b,w}$-(semi)stable if and only if
	\begin{itemize}
		\item $E[k]\in\coh^b(X)$ for some $k\in\Z$, and
		\item $\nu\_{b,w}(F)\, < (\leq)\,\nu\_{b,w}\big(E[k]/F\big)$ for all non-trivial subobjects $F\hookrightarrow E[k]$ in $\coh^b(X)$.
	\end{itemize}
\end{Def}
\noindent Note that $\nu\_{b,w}$-stability defines a Harder-Narasimhan (HN) filtration on $\Coh^b(X)$, and so a \emph{weak stability condition} on $\cD(X)$ \cite[Proposition 12.2]{bayer:the-space-of-stability-conditions-on-abelian-threefolds}. 

In Figure \ref{projetcion} we plot the $(b,w)$-plane simultaneously with the image of the projection map
\begin{eqnarray}\label{pi}
	\Pi\colon\ K(X) \setminus \big\{E \colon \ch_0(E) = 0\big\}\! &\longrightarrow& \mathbb{R}^2, \\
	E &\ensuremath{\shortmid\joinrel\relbar\joinrel\rightarrow}& \!\!\bigg(\frac{\ch_1(E).H^2}{\ch_0(E)H^3}\,,\, \frac{\ch_2(E).H}{\ch_0(E)H^3}\bigg). \nonumber
\end{eqnarray}
\begin{figure}[h]
	\begin{centering}
		\definecolor{zzttqq}{rgb}{0.27,0.27,0.27}
		\definecolor{qqqqff}{rgb}{0.33,0.33,0.33}
		\definecolor{uququq}{rgb}{0.25,0.25,0.25}
		\definecolor{xdxdff}{rgb}{0.66,0.66,0.66}
		
		\begin{tikzpicture}[line cap=round,line join=round,>=triangle 45,x=1.0cm,y=1.0cm]
		
		\draw[->,color=black] (-4,0) -- (4,0);
		\draw  (4, 0) node [right ] {$b,\, \frac{H^{n-1}\ch_1}{H^n\ch_0}$};


		\fill [fill=gray!40!white] (0,0) parabola (3,4) parabola [bend at end] (-3,4) parabola [bend at end] (0,0);
		
		\draw  (0,0) parabola (3.1,4.27); 
		\draw  (0,0) parabola (-3.1,4.27); 
		\draw  (3.8 , 3.6) node [above] {$w= \frac{b^2}{2}$};
		
		

		\draw[->,color=black] (0,-.8) -- (0,4.7);
		\draw  (1, 4.3) node [above ] {$w,\,\frac{H^{n-2}\ch_2}{H^n\ch_0}$};

		
		\draw [dashed, color=black] (2.3,1.5) -- (2.3,0);
		\draw [dashed, color=black] (2.3, 1.5) -- (0, 1.5);
		\draw [color=black] (2.6, 1.36) -- (-1.3, 3.14);
		
		\draw  (2.8, 1.8) node {$\Pi(E)$};
		\draw  (1, 3) node [above] {\Large{$U$}};
		\draw  (0, 1.5) node [left] {$\frac{H^{n-2}\ch_2}{H^n\ch_0}$};
		\draw  (2.3 , 0) node [below] {$\frac{H^{n-1}\ch_1}{H^n\ch_0}$};
		\begin{scriptsize}
		\fill (0, 1.5) circle (1.4pt);
		\fill (2.3,0) circle (1.4pt);
		\fill (2.3,1.5) circle (1.4pt);
		\fill (-1,3) circle (1.4pt);
		\draw  (-1.2, 2.96) node [below] {$(b,w)$};
		
		\end{scriptsize}
		
		\end{tikzpicture}
		
		\caption{$(b,w)$-plane and the projection $\Pi(E)$ when $\ch_0(E)>0$}
		
		\label{projetcion}
		
	\end{centering}
\end{figure}

\noindent The pair $(b,w)$ for our weak stability condition $\nu\_{b,w}$ is in the shaded open subset
\begin{equation}\label{Udef}
U \,:= \,\left\{(b,w) \in \mathbb{R}^2 \colon w > \tfrac12b^2  \right\}.
\end{equation}
Also, the image $\Pi(E)$ of $\nu_{b,w}$-semistable objects $E$ with $\ch_0(E)\ne0$ is always \emph{outside} $U$ due to the classical Bogomolov-Gieseker-type inequality
for the $H$-discriminant,
\begin{equation}\label{discr}
\Delta_H(E)\ =\  \big(\!\ch_1(E).H^2\big)^2 -2 (\ch_0(E)H^3)(\ch_2(E).H)\ \ge\ 0,
\end{equation}
proved for $\nu\_{b,w}$-semistable objects $E$ in \cite[Theorem 3.5]{bayer:the-space-of-stability-conditions-on-abelian-threefolds}. 
The slope $\nu_{b,w}(E)$ of $E$ is the gradient of the line connecting $(b,w)$ to $\Pi(E)$.

Objects of $\cD(X)$ give the space of weak stability conditions a wall and chamber structure, explained in \cite[Proposition 4.1]{feyz:noether-lefschtz-loci} for instance.

\begin{Prop}[\textbf{Wall and chamber structure}]\label{prop. locally finite set of walls}
	Fix an object $E \in \mathcal{D}(X)$ such that the vector $\left(\ch_0(E), \ch_1(E).H^2, \ch_2(E).H\right)$ is non-zero. There exists a set of lines $\{\ell_i\}_{i \in I}$ in $\mathbb{R}^2$ such that the segments $\ell_i\cap U$ (called ``\emph{walls}") are locally finite and satisfy 
	\begin{itemize*}
		\item[\emph{(}a\emph{)}] If $\ch_0(E)\ne0$ then all lines $\ell_i$ pass through $\Pi(E)$.
		\item[\emph{(}b\emph{)}] If $\ch_0(E)=0$ then all lines $\ell_i$ are parallel of slope $\frac{\ch_2(E).H^{n-2}}{\ch_1(E).H^{n-1}}$.
		\item[\emph{(}c\emph{)}] The $\nu\_{b,w}$-(semi)stability of $E$ is unchanged as $(b,w)$ varies within any connected component (called a ``\emph{chamber}") of $U \setminus \bigcup_{i \in I}\ell_i$.
		\item[\emph{(}d\emph{)}] For any wall $\ell_i\cap U$ there is $k_i \in \mathbb{Z}$ and a map $f\colon E[k_i] \to F$ in $\cD(X)$ such that
		\begin{itemize}
			\item for any $(b,w) \in \ell_i \cap U$, the objects $E[k_i],\,F$ lie in the heart $\cA(b)$,
			\item $E[k_i]$ is $\nu\_{b,w}$-semistable with $\nu\_{b,w}(E)=\nu\_{b,w}(F)=\,\mathrm{slope}\,(\ell_i)$ constant on the wall $\ell_i \cap U$, and
			\item $f$ is a surjection $E[k_i] \twoheadrightarrow F$ in $\cA(b)$ which strictly destabilises $E[k_i]$ for $(b,w)$ in one of the two chambers adjacent to the wall $\ell_i$.
			\hfill$\square$
			
		\end{itemize} 
	\end{itemize*} 
\end{Prop}

\begin{figure} [h]
	\begin{centering}
		
		\begin{tikzpicture}[line cap=round,line join=round,>=triangle 45,x=1.0cm,y=1.0cm]
		
		\draw[->,color=black] (-10.5,0) -- (-5.5,0);
		\draw[->,color=black] (-3,0) -- (2,0);
		
		\fill [fill=gray!30!white] (-0.5,0) parabola (1.47, 3.03) parabola [bend at end] (-2.47,3.03) parabola [bend at end] (-0.5,0);
		
		\fill [fill=gray!30!white] (-8,0) parabola (-6.03, 3.03) parabola [bend at end] (-9.97,3.03) parabola [bend at end] (-8,0);

		\draw[->,color=black] (-8,-1) -- (-8,3.5);
		\draw[->,color=black] (-0.5,-1) -- (-0.5,3.5);

		\draw [] (-0.5,0) parabola (1.5,3.12); 
		\draw [] (-0.5,0) parabola (-2.5,3.12); 
		\draw [] (-8,0) parabola (-10,3.12); 
		\draw [] (-8,0) parabola (-6,3.12);

		\draw[color=black, dashed] (-10.5,2.8) -- (-6,1);
		\draw[color=black, dashed] (-10.5,1.8) -- (-6.5,0.2);

		\draw[color=black,semithick] (-9.8,2.52) -- (-6.7,1.28);
		\draw[color=black,semithick] (-9.3,1.32) -- (-7.2,.48);
		
		\draw (-10.5,1.8) node [left] {$\ell_2$};
		\draw (-10.5,2.8) node [left] {$\ell_1$};

		\draw (.8,3.5) node [right] {\large{$\ch_0(E) \neq 0$}};
		\draw (-6.8,3.5) node [right] {\large{$\ch_0(E) = 0$}};
		\draw (-5.5,0) node [below] {$b, \frac{H^{n-1}\ch_1}{H^n\ch_0}$};
		\draw (-8,3.5) node [above] {$w, \frac{H^{n-2}\ch_2}{H^n\ch_0}$};
		
		\draw (2,0) node [below] {$b, \frac{H^{n-1}\ch_1}{H^n\ch_0}$};
		\draw (-0.5,3.5) node [above] {$w, \frac{H^{n-2}\ch_2}{H^n\ch_0}$};
		
		\draw (-7.3,2.5) node [right] {\Large{$U$}};
		\draw (0.2,2.5) node [right] {\Large{$U$}};

		\draw (1.5, 1.2) node [right] {$\Pi(E)$};
		
		\draw[color=black, dashed] (1.5, 1.2) -- (-2.9,2.8);
		\draw[color=black, dashed] (1.5, 1.2) -- (-2.2, .7);
		
		\draw (-2.9,2.8) node [left] {$\ell_1$};
		\draw (-2.2, .7) node [left] {$\ell_2$};

		\draw[color=black, semithick] (-2.3 ,2.58) -- (.88,1.423);
		\draw[color=black, semithick] (-1.5,.795) -- (0.7,1.092);

		\begin{scriptsize}
		
		%
		%
		%

		\fill [color=black] (1.5,1.2) circle (1.4pt);
		
		%

		
		%
		%
		%
		%
		%
		%
		

		\end{scriptsize}
		
		\end{tikzpicture}
		
		\caption{The line segments $\ell_i \cap U$ are walls for $E$.}
		
		\label{wall.figure}
		
	\end{centering}
	
\end{figure}

\section{large volume limit}

As before $X$ is a smooth projective complex variety of dimension $n \geq 2$ and $H$ is an ample divisor on $X$. 
In this section, we first prove the following general statement which holds for any coherent sheaf. 
	\begin{Thm}\label{thm.bounding giesker chamber}
	Given a coherent sheaf $E$ of rank bigger than one. There is no wall for $E$ crossing the vertical line $b = b_0$ at a point inside $U$ whenever 
	\begin{equation}\label{interval of b}
	\mu(E) - \frac{1}{H^n\ch_0(E) (\ch_0(E) -1)} \leq b_0 < \mu(E) \ . 
	\end{equation} 
\end{Thm}

\begin{proof}
	Suppose for a contradiction that there is a wall $\ell$ for $E$ which intersects the vertical line $b=b_0$. Since $E$ is $\nu_{b,w}$-semistable for $(b,w) \in \ell \cap U$, we have $\Delta_H(E) \geq 0$ \eqref{discr}. If $\Delta_H(E) = 0$, then \cite[Corollary 3.11(a)]{bayer:the-space-of-stability-conditions-on-abelian-threefolds} implies that there is no wall for $E$ in $U$, so we may assume $\Delta_H(E) > 0$. Hence $\Pi(E)$ lies outside $\overline{U} = U \cup \partial U$. Let
	\begin{equation*}
	(r, c, h) \coloneqq  \left(H^n\ch_0(E),\, H^{n-1}\ch_1(E),\, H^{n-2}\ch_2(E)  \right),
	\end{equation*} 
    and let $\delta_i = \frac{c}{r} - b_i$ where $b_2 < b_1$ are the $b$-values at the intersection points of the wall $\ell$ with $\partial U$, see Figure \ref{large volume limit.shape}.  
		\begin{figure} [h]
		\begin{centering}
			
			\begin{tikzpicture}[line cap=round,line join=round,>=triangle 45,x=1.0cm,y=1.0cm]
			
			
			\fill [fill=gray!30!white] (0,0) parabola (2,4) parabola [bend at end] (-2,4) parabola [bend at end] (0,0);
			
			\draw (3,0) node [right] {$b, \frac{H^{n-1}\ch_1}{H^n\ch_0}$};
			\draw (0,4.2) node [above] {$w, \frac{H^{n-2}\ch_2}{H^n\ch_0}$};
			
			\draw[->,color=black] (-3,0) -- (3,0);

			\draw [] (0,0) parabola (2,4); 
			\draw [] (0,0) parabola (-2,4);

			\draw [ color=black, dashed] (1.4,0)--(1.4,1.95);
			\draw [ color=black, dashed] (1.9,1.9)--(1.9,0);
			\draw [dashed, color= black] (1, 4)--(1,0);
			\draw [dashed, color= black] (-1.5,2.25)--(-1.5,0);
			
			\draw [dashed, color=black] (1.9,1.9) -- (1.4,1.95);
			\draw [thick] (1.4,1.95) -- (-1.5,2.25);

			\draw[->] (0,0) -- (0,4.2);
			
			\draw (1.49,0) node [below] {$b_1$};
			\draw (-1.5,0) node [below] {$b_2$};
			\draw (1,0) node [below] {$b_0$};
			\draw (-.5,2.1) node [above] {$\ell$};
			\draw (2.55,1.6) node [above] {$\Pi(E)$};
			\draw (0,0) node [below] {$o$};
			
			\draw  (1.9,0) node [below ] {$\frac{c}{r}$};

			\begin{scriptsize}
			\fill (0, 0) circle (1.3pt);
			
			\fill (1.9,0) circle (1.3pt);
			\fill [color=black] (1.9, 1.9) circle (1.3pt);
			
			\fill [color=black] (-1.5,2.25) circle (1.3pt);
			\fill [color=black] (1.4,1.95) circle (1.3pt);
				\fill [color=black] (1, 0) circle (1.3pt);
			\fill [color=black] (1.4,0) circle (1.3pt);
			\fill [color=black] (-1.5,0) circle (1.3pt);
			
			\end{scriptsize}
			
			\end{tikzpicture}
			
			\caption{The wall $\ell$}
			
			\label{large volume limit.shape}
			
		\end{centering}
		
	\end{figure}

	\noindent  The inequality \eqref{interval of b} gives
	\begin{equation}\label{assumption-1}
	\frac{c}{r} - \frac{H^n}{r(r-H^n)} \leq b_0 < b_1   
	\end{equation}
	which implies $0 < \delta_1 = \frac{c}{r} -b_1 < \frac{H^n}{r(r-H^n)}$. 
	 One can easily check  
	\begin{equation*}
	\delta_1 \delta_2 = \frac{\Delta_H(E)}{r^2}. 
	\end{equation*}
	Thus we get
	\begin{align}\label{b-diff}
	b_1 -b_2 = \delta_2 -\delta_1  & =   \frac{\Delta_H(E)}{ \delta_1 r^2}- \delta_1 \nonumber\\
	& > \frac{\Delta_H(E)}{r^2}. \frac{r(r-H^n)}{H^n} - \frac{H^n}{r(r-H^n)}\ .  
	\end{align}
	\textbf{Step 1.} First assume $r \geq 3H^n$ or $\Delta_H(E) \geq 3$, then \eqref{b-diff} gives
	\begin{equation*}
	b_1 -b_2 > \frac{1}{r-H^n}\ .
	\end{equation*}    	
	Thus there is $k \in \mathbb{Z}$ such that 
	\begin{equation*}
	b_2 \ < b^* \coloneqq  \frac{k}{r-H^n} < \ b_1\ ,  
	\end{equation*}
	so the wall $\ell$ intersects the vertical line $b= b^*$ at a point $(b^*, w^*)$ inside $U$. Let 
	\begin{equation}\label{sequence}
	E_1 \hookrightarrow E \twoheadrightarrow E_2
	\end{equation}
	be a destabilising sequence for $E$ which makes it unstable in one side of the wall $\ell$ with 
	\begin{equation*}
	(r_i, c_i, h_i) \coloneqq \left(H^n\ch_0(E_i),H^{n-1}\ch_1(E_i), H^{n-2}\ch_2(E_i)  \right). 
	\end{equation*}
At any $(b,w) \in \ell \cap U$, we have $\nu\_{b,w}(E) = \nu\_{b,w}(E_i)$ for $i=1,2$, so 
\begin{equation}\label{imaginary part}
0 < \text{Im}[Z_{b,w}(E_i)] = c_i -b\,r_i < \text{Im}[Z_{b,w}(E)] = c-b\,r. 
\end{equation}   
We consider two different cases:

	\textbf{Case i.} First assume 
	\begin{equation*}
	\frac{h_1-r_1  w}{c_1 -r_1 b} = \nu\_{b,w}(E_1) < \nu\_{b,w}(E) = \frac{ \tfrac{h}{r} - w}{\tfrac{c}{r} -b } \ .
	\end{equation*} 
	for $(b,w) \in U$ above the wall $\ell_E$, i.e. our destabilising sequence destabilises $E$ below $\ell$. Then $r_1 \neq 0$ and
	\begin{equation*}
	\lim\limits_{w \rightarrow +\infty} \nu\_{b,w}(E_1) = \frac{-r_1}{c_1 -br_1}\  <\  \lim\limits_{w \rightarrow +\infty} \nu\_{b,w}(E) = \frac{-1}{\frac{c}{r} -b}. 
	\end{equation*}
	Suppose $b_2 < b <  b_1$, then \eqref{imaginary part} implies $r_1 > 0$ and so
	\begin{equation}\label{ordering of slope}
	\frac{c_1}{r_1} < \frac{c}{r}\ . 
	\end{equation} 
	The limit of \eqref{imaginary part} when $b \rightarrow b_1^{-}$ implies that $0 \leq c_1 - r_1 b_1$.    
	Combining this with \eqref{ordering of slope} gives 
	\begin{equation}\label{u-l-b}
	0 < \frac{c}{r} - \frac{c_1}{r_1} \ \leq\ \delta_1 = \frac{c}{r} -b_1\ .  
	\end{equation}
	If $0 < r_1 < r$, then 
	\begin{equation}\label{bounding.1}
	\dfrac{c}{r} - \dfrac{c_1}{r_1} = \dfrac{cr_1-c_1r}{rr_1} = \dfrac{cH^n\text{ch}_0(E_1) - c_1H^n\text{ch}_0(E) }{rr_1} \geq \dfrac{H^n}{r(r-H^n)}\ , 
	\end{equation}  
	and if $r \leq r_1$, inequality \eqref{imaginary part} at $b= b^*$ gives
	\begin{align}\label{bounding.2}
	\frac{c}{r} - \frac{c_1}{r_1} = \frac{c-b^*r}{r} - \frac{c_1-b^*r_1}{r_1} \geq \frac{c-b^*r}{r} - \frac{c_1-b^*r_1}{r} \geq \frac{1}{r} \left( \dfrac{H^n}{r-H^n} \right).     
	\end{align}
	The last inequality comes from the point that $b^*r$ and $b^*r_1$ are integral multiples of  $\frac{1}{\ch_0(E) -1}$. 
Therefore, by \eqref{u-l-b} we get
	\begin{equation*}
	\delta_1 = \frac{c}{r} -b_1 \geq \frac{H^n}{r(r-H^n)}
	\end{equation*}  
	which is not possible by our assumption \eqref{assumption-1}.  
	
	\textbf{Case ii.} Suppose our destabilising sequence makes $E$ unstable above the line $\ell_E$, then    
	\begin{equation*}
	\frac{h_2-r_2  w}{c_2-r_2 b} = \nu\_{b,w}(E_2) < \nu\_{b,w}(E) = \frac{ \tfrac{h}{r} - w}{\tfrac{c}{r} -b } 
	\end{equation*} 
	for $(b,w) \in U$ above $\ell_E$. Hence $r_2 \neq 0$ and
	\begin{equation*}
	\lim\limits_{w \rightarrow +\infty} \nu\_{b,w}(E_2) = \frac{-r_2}{c_2 -r_2b}\  <\  \lim\limits_{w \rightarrow +\infty} \nu\_{b^*,w}(E) = \frac{-1}{\frac{c}{r} -b^*}
	\end{equation*}
	If $b_2<b <b_1$, \eqref{imaginary part} implies $r_2 >0$ and so
	\begin{equation}
	\frac{c_2}{r_2} < \frac{c}{r}\ . 
	\end{equation}
	Then limit of \eqref{imaginary part} when $b \rightarrow b_1^-$ implies that 
	\begin{equation*}
	\delta_1 = \frac{c}{r} -b_1 \geq \frac{c}{r} - \frac{c_2}{r_1} > 0\ .
	\end{equation*}
	Hence applying the same argument as in \eqref{bounding.1} and \eqref{bounding.2} reach a contradiction.
	
	\textbf{Step 2.} Finally, we show that if $r=2H^n$ and $\Delta_H(E) = 1$, or $2$, there is no wall for $E$ in $U$, so the claim follows. Suppose there is a wall with the destabilising sequence \eqref{sequence}. The same argument as in \cite[Corollary 3.10]{bayer:the-space-of-stability-conditions-on-abelian-threefolds} gives 
	\begin{equation*}
	0 \leq \Delta_H(E_i) = c_i^2 - 2r_ih_i < \Delta_H(E)\ . 
	\end{equation*}
	If $\Delta_H(E) =1$, then $\Delta_H(E_i) =0$ for $i=1, 2$, so
	\begin{align*}
	1 = \Delta_H(E) = (c_1+c_2)^2 - 2(r_1+r_2)(h_1+h_2) = 2c_1c_2 -2(r_1h_2 +r_2h_1) 
	\end{align*}
	which is not possible. Similarly, if $\Delta_H(E) = 2$, one of the following cases happens: 
	\begin{enumerate*}
		\item $\Delta_H(E_i) = 0$ for $i=1, 2$. Since $H^n|r_i$, we get $H^n|c_i^2$, but 
		\begin{align*}
		1 = \frac{\Delta_H(E)}{2} = c_1c_2 - (r_1h_2 +r_2h_1) 
		\end{align*}
		implies that $\gcd(H^n, c_i) = 1$, a contradiction.
		\item  $\Delta_H(E_1) = 1$ and $\Delta_H(E_2) = 0$, then
		\begin{align*}
		1 = \Delta_H(E)- \Delta_H(E_1) =  2c_1c_2 -2(r_1h_2 +r_2h_1) 
		\end{align*}
		which is not again possible. Similar argument applies to the case $\Delta_H(E_1) = 0$ and $\Delta_H(E_2) = 1$. 
		\item $\Delta_H(E_i) =1$ for $i=1, 2$, then $\gcd(H^n, c_i^2) =1$ for $i=1, 2$, but 
		\begin{equation*}
		0 = \Delta_H(E) - \Delta_H(E_1) -\Delta_H(E_2) = 2c_1c_2 - 2(r_1h_2 +r_2h_1) 
		\end{equation*}  
		which implies $H^n | c_1c_2$, a contradiction. 
	\end{enumerate*} 
\end{proof}	
By combing the ideas in the proof of Theorem \ref{thm.bounding giesker chamber} and the notion of safe area introduced in \cite{feyz:rank-r-dt-theory}, we get the following. 
\begin{Prop}\label{prop.line of stability-E}
	Take a $\mu$-stable coherent sheaf $E$ of positive rank, 
	and let $\ell$ be a line passing through $\Pi(E)$ which intersects $\partial U$ at two points with $b$-values $b_2 <b_1 \leq \mu(E)$ satisfying
	\begin{equation*}
  0\leq \mu(E) - b_1 < \mu(E)-\mu^{\max}(E) \qquad \text{and} \qquad \ch_0(E)(\mu(E) - b_1)< b_1-b_2\,, 
	\end{equation*}
	where $\mu^{\max}(E)$ is defined as in \eqref{max}. Then $E$ is $\nu\_{b,w}$-stable for any $(b,w) \in \ell \cap U$.  
\end{Prop}
\begin{proof}
A similar argument as in \cite[Proposition 14.2]{bridgeland:K3-surfaces} or \cite[Lemma 2.7]{bayer:the-space-of-stability-conditions-on-abelian-threefolds} implies that $E$ is $\nu_{b,w}$-stable when $b< \mu(E)$ and $w \gg 0$. 
	
If $\Delta_H(E) = 0$, then \cite[Corollary 3.11(a)]{bayer:the-space-of-stability-conditions-on-abelian-threefolds} implies that there is no wall for $E$ in $U$, so the claim follows. Thus we may assume $\Delta_H(E) > 0$, i.e. $\Pi(E)$ lies outside $\overline{U} = U \cup \partial U$. Then there is a unique line through $\Pi(E)$ which is tangent to $\partial U$ at a point to the left of $\Pi(E)$. If we move this line upwards while making it pass through $\Pi(E)$, its intersection point with $\partial U$ move further apart until we find a unique $\ell_E$ for which the $b$-values $b_2^* < b_1^*$ of the two points of $\ell_E \cap \partial U$ satisfy 
\begin{equation}\label{def-ell}
\mu(E) -b_1^* = \min \left\{  \mu(E)-\mu^{\max}(E)\ \ ,\ \ \frac{1}{\ch_0(E)}(b_1^* -b_2^*)  \right\}. 
\end{equation}
We show there is no wall for $E$ above the line $\ell_E$ to the left of $\Pi(E)$.  Since the line $\ell$ in the Proposition lies above $\ell_E$, the claim follows.
		\begin{figure} [h]
	\begin{centering}
		
		\begin{tikzpicture}[line cap=round,line join=round,>=triangle 45,x=1.0cm,y=1.0cm]
		
		
		\fill [fill=gray!30!white] (0,0) parabola (2,4) parabola [bend at end] (-2,4) parabola [bend at end] (0,0);
		
		\draw (3,0) node [right] {$b, \frac{H^{n-1}\ch_1}{H^n\ch_0}$};
		\draw (0,4.2) node [above] {$w, \frac{H^{n-2}\ch_2}{H^n\ch_0}$};
		
		\draw[->,color=black] (-3,0) -- (3,0);

		\draw [] (0,0) parabola (2,4); 
		\draw [] (0,0) parabola (-2,4);

		\draw [ color=black, dashed] (1.4,0)--(1.4,1.95);
		\draw [ color=black, dashed] (1.9,1.9)--(1.9,0);
		\draw [dashed, color= black] (-1.5,2.25)--(-1.5,0);
		
		\draw [dashed, color=black] (1.9,1.9) -- (1.4,1.95);
		\draw [thick] (1.4,1.95) -- (-1.5,2.25);
		
		\draw [dashed] (1.9,1.9) -- (1.46,2.11);
		\draw [color=black, thick]  (1.46,2.11) -- (-1.9,3.7);

		\draw[->] (0,0) -- (0,4.2);
		
		\draw (1.45,0) node [below] {$b_1^*$};
		\draw (-1.5,0) node [below] {$b_2^*$};
		\draw (-.5,2.1) node [above] {$\ell_E$};
		\draw (-.5,3.) node [above] {$\tilde{\ell}$};
		\draw (2.55,1.6) node [above] {$\Pi(E)$};
		\draw (0,0) node [below] {$o$};
		
		\draw  (2.1,.05) node [below ] {\small $\mu(E)$};

		\begin{scriptsize}
		\fill (0, 0) circle (1.3pt);
		
		\fill (1.9,0) circle (1.3pt);
		\fill [color=black] (1.9, 1.9) circle (1.3pt);
		
		\fill [color=black] (-1.5,2.25) circle (1.3pt);
		\fill [color=black] (1.4,1.95) circle (1.3pt);
		\fill [color=black] (1.4,0) circle (1.3pt);
		\fill [color=black] (-1.5,0) circle (1.3pt);
		\fill [color=black] (-1.92,3.7) circle (1.3pt);
     	\fill [color=black] (1.46,2.11) circle (1.3pt);
		
		\end{scriptsize}
		
		\end{tikzpicture}
		
		\caption{The line $\ell_E$}
		
		\label{fig.line-ell-E}
		
	\end{centering}
	
\end{figure}

Suppose for a contradiction there is a wall $\tilde{\ell}$ for $E$ which lies above $\ell_E$ to the left of $\Pi(E)$ and intersects $\partial U$ at two points with $b$-values $\tilde{b}_2< \tilde{b}_1$ satisfying
\begin{equation}\label{inter}
\tilde{b}_2< b_2^*    \qquad \text{and} \qquad b_1^*< \tilde{b}_1 < \mu(E). 
\end{equation}
Let $E_1 \hookrightarrow E \twoheadrightarrow E_2$ be a destabilising sequence for $E$ which make it unstable below $\tilde{\ell}$.  We first show that both $E_1$ and $E_2$ are coherent sheaves. Taking cohomology from the destabilising sequence gives the long exact sequence of coherent sheaves   	
\begin{equation*}
0 \rightarrow \cH^{-1}(E_2) \rightarrow \cH^0(E_1) \rightarrow \cH^0(E) \rightarrow \cH^0(E_2) \rightarrow 0 \ . 
\end{equation*}	
This immediately gives $E_1$ is a sheaf. Suppose $\cH^{-1}(E_2)$ is of rank $r' \neq 0$. As we move $b \uparrow \tilde{b}_1$ or $b \downarrow \tilde{b}_2$ along $\tilde{\ell}$, $E_1$ and $E_2$ remain in the heart $\Coh^b(X)$ by Proposition \ref{prop. locally finite set of walls}. This implies
\begin{equation*}
\mu(\cH^{-1}(E_2)) \leq \tilde{b}_2  \qquad \text{and} \qquad \mu(\cH^0(E_i)) \geq \tilde{b}_1 \,\,\,\,\,\,\ \text{for $i=1, 2$}.
\end{equation*} 
Combining this with \eqref{inter} gives 
\begin{align*}
\ch_1(E)H^{n-1} = & \ \ch_1(E_1)H^{n-1} +\ch_1(\cH^{0}(E_2))H^{n-1} -\ch_1(\cH^{-1}(E_2))H^{n-1} \\
\geq & \ H^n\tilde{b}_1 (\ch_0(E) + r') - H^n\tilde{b}_2 r'\\
> & \ H^n\,b_1^*\ch_0(E) +  H^nr' (b_1^* -b_2^*).
\end{align*}	
Therefore $\mu(E) -b_1^* > \frac{1}{\ch_0(E)} (b_1^* -b_2^*)$ which is not possible by \eqref{def-ell}, so $r' =0$. We know $\cH^{-1}(E_2)$ is a torsion-free sheaf by the definition \eqref{Abdef} of $\coh^b(X)$, thus it is zero and $E_2$ is a sheaf. Hence $E_1$ is a subsheaf of $E$. The wall $\tilde{\ell}$ intersects the vertical line $b = b_1^*$ at a point $(b^*, w)$ inside $U$ where $E_1$ and $E$ have the same $\nu_{b^*_1, w}$-slope, see Figure \ref{fig.line-ell-E}. Thus 
\begin{equation*}
0 < \text{Im}[Z_{b^*_1, w}(E_1)] = \ch_1(E_1)H^{n-1}-b_1^*\ch_0(E_1)H^n 
\end{equation*} 
which implies 
\begin{equation*}
 b_1^* < \mu(E_1) < \mu(E).
\end{equation*}
The second inequality comes from $\mu$-stability of $E$. Note that the quotient sheaf $E_2$ cannot be supported in co-dimension at least 2, otherwise its $\nu_{b,w}$-slope is $+\infty$ and it cannot destabilise $E$, so $\mu(E_1) \neq \mu(E)$. Therefore
\begin{equation*}
\mu(E) -b_1^* > \mu(E) -\mu(E_1) \geq \mu(E) -\mu^{\max}(E)
\end{equation*} 
 which is not possible by \eqref{def-ell}.  	
\end{proof}

The next step is to analyse points $(b,w) \in U$ where the shift of a slope-stable reflexive sheaf is stable. The following Lemma is well-known, we add it for completeness.

\begin{Lem}\label{lem. phase.1.locally free sheaaf}
	Let $E$ be a $\mu$-stable reflexive sheaf of positive rank. Then $E[1]$ is $\nu_{b_0,w}$-stable for $b_0=\mu(E)$ and all $w > \frac{b_0^2}{2}$.     
\end{Lem}
\begin{proof}
	By definition $E[1] \in \coh^{b_0}(X)$ and Im$[Z_{b_0, w}(E)] = 0$. Therefore $E[1]$ is $\nu_{b_0, w}$-semistable of slope $+ \infty$. Suppose $E[1]$ is strictly $\nu_{b_0,w}$-semistable, so there is a short exact sequence in $\coh^{b_0}(X)$
	\begin{equation}\label{exact.1.}
	E_1 \hookrightarrow E[1] \twoheadrightarrow E_2
	\end{equation}
	such that $\nu\_{b_0, w}(E_1) = \nu\_{b_0, w}(E_2) = +\infty$. 
	
	\textbf{Case i.} First assume $\cH^{-1}(E_1) = 0$, then 
	\begin{equation}\label{im}
	\text{Im}[Z\_{b_0, w}(E_1)] = \ch_1(\cH^{0}(E_1)).H^{n-1}-b_0 H^n\ch_0(\cH^{0}(E_1)) =0\ . 
	\end{equation}	
	If $\ch_0(\cH^{0}(E_1)) \neq 0$, we get $\mu(\cH^0(E_1)) =b_0$ which is not possible by the definition \eqref{Abdef} of $\coh^{b_0}(X)$. Therefore, $E_1$ is a sheaf supported in co-dimension at least 2. Taking cohomology from \eqref{exact.1.} gives the short exact sequence of coherent sheaves
	\begin{equation}\label{non-trivila ext}
	0 \rightarrow E \rightarrow \cH^{-1}(E_2) \rightarrow \cH^{0}(E_1) \rightarrow 0\ .
	\end{equation} 
	Since $E$ is a reflexive sheaf, \cite[Proposition 1.1]{hartshorne:stable-reflexive-sheaves} implies that $E$ lies in the short exact sequence
	\begin{equation*}
	0 \rightarrow E \rightarrow G \rightarrow G' \rightarrow 0 
	\end{equation*}
	where $G$ is a locally free sheaf and $G'$ is torsion-free. Since $\cH^0(E_1)$ is supported in co-dimension at least 2, we have
	\begin{equation*}
 \text{Hom}(\cH^0(E_1), G') = 0 = \text{Ext}^1(\cH^0(E_1), G). 
	\end{equation*}
	This implies Ext$^1(\cH^0(E_1), E) = 0$, thus \eqref{non-trivila ext} gives $\cH^{-1}(E_2) \cong E \oplus \cH^0(E_1)$ which is not possible because $\cH^{-1}(E_2)$ is torsion-free.

	\textbf{Case ii.} If $\cH^{-1}(E_1) \neq 0$, the exact sequence $\cH^{-1}(E_1)[1] \hookrightarrow E_1 \twoheadrightarrow \cH^0(E_1)$ 
	in $\coh^{b_0}(X)$ implies that $\cH^{-1}(E_1)[1]$ is a subobject of $E[1]$ of $\nu\_{b_0, w}$-slope $+\infty$; so there is a short exact sequence
	\begin{equation}\label{s.1}
	\cH^{-1}(E_1)[1] \hookrightarrow E[1] \twoheadrightarrow E_2'
	\end{equation}  
	in $\coh^{b_0}(X)$ where $E_2'$ lies in the exact sequence $\cH^0(E_1) \hookrightarrow E_2' \twoheadrightarrow E_2$ in $\Coh^{b_0}(X)$. Note that $E_2' \neq 0$ otherwise $E_2 \cong \cH^0(E_1)[1]$ which is not possible because $\cH^0(E_1)[1] \notin \Coh^{b_0}(X)$. 
	
	Taking cohomology from \eqref{s.1} gives a short exact sequence of coherent sheaves  
	\begin{equation*}
	0 \rightarrow \cH^{-1}(E_1) \rightarrow  E \rightarrow  \cH^{-1}(E_2') \rightarrow 0,
	\end{equation*} 
	thus $\cH^0(E_2') = 0$. The torsion-free sheaves $E$, $\cH^{-1}(E_1)$ and $\cH^{-1}(E_2')$ satisfy 
	\begin{equation*}
	\text{Im}\big[Z\_{b_0, w}\big(\cH^{-1}(E'_2)\big)\big] = 	\text{Im}\big[Z\_{b_0, w}\big(\cH^{-1}(E_1)\big)\big] = 	\text{Im}\big[Z\_{b_0, w}\big(\cH^{-1}(E)\big)\big] = 0. 
	\end{equation*}
	Thus they all have the same $\mu$-slope which is not possible by $\mu$-stability of $E$. 
\end{proof}
By applying a similar argument to the proof of Proposition \ref{prop.line of stability-E}, we get the following.  
\begin{Prop}\label{prop.line of stability-E[1]}
	Take a $\mu$-stable reflexive sheaf $E$ of positive rank, and let 
	$\ell$ be a line passing through $\Pi(E)$ which intersects $\partial U$ at two points with $b$-values $b_2 <b_1 $ such that
	\begin{equation*}
	0 \leq b_2 - \mu(E) < \mu^{\min}(E)-\mu(E) \qquad \text{and} \qquad \ch_0(E)(b_2- \mu(E))< b_1-b_2\,, 
	\end{equation*}
	where $\mu^{\min}(E)$ is defined as in \eqref{min}. Then $E[1]$ is $\nu_{b,w}$-stable for any $(b,w) \in \ell \cap U$. 
\end{Prop}
\begin{proof}
	Lemma \ref{lem. phase.1.locally free sheaaf} implies $E[1]$ is $\nu_{b_0, w}$-stable for $b_0 = \mu(E)$. Thus the structure of walls described in Proposition \ref{prop. locally finite set of walls} implies that $E[1]$ is $\nu_{b, w}$-stable where $\mu(E)<b$ and $w \gg 0$.
	
	If $\Delta_H(E) = 0$, there is no wall for $E[1]$ in $U$ by \cite[Corollary 3.11(a)]{bayer:the-space-of-stability-conditions-on-abelian-threefolds}, so the claim follows. Hence we may assume  $\Delta_H(E) > 0$. There is a unique line $\ell_{E[1]}$ passing through $\Pi(E)$ such that the $b$-values $b_2^*< b_1^*$ of the two points of $\ell_{E[1]} \cap \partial U$ satisfy 
     \begin{equation}\label{def-ell-[1]}
     b_2^*-\mu(E) = \min \left\{  \mu^{\min}(E)-\mu(E)\ \ ,\ \ \frac{1}{\ch_0(E)}(b_1^* -b_2^*)  \right\}.
     \end{equation}
     		\begin{figure} [h]
     	\begin{centering}
     		
     		\begin{tikzpicture}[line cap=round,line join=round,>=triangle 45,x=1.0cm,y=1.0cm]
     		
     		
     		\fill [fill=gray!30!white] (0,0) parabola (2,4) parabola [bend at end] (-2,4) parabola [bend at end] (0,0);
     		
     		\draw (3,0) node [right] {$b, \frac{H^{n-1}\ch_1}{H^n\ch_0}$};
     		\draw (0,4.2) node [above] {$w, \frac{H^{n-2}\ch_2}{H^n\ch_0}$};
     		
     		\draw[->,color=black] (-3,0) -- (3,0);

     		\draw [] (0,0) parabola (2,4); 
     		\draw [] (0,0) parabola (-2,4);

     		\draw [ color=black, dashed] (-1.4,0)--(-1.4,1.95);
     		\draw [ color=black, dashed] (-1.9,1.9)--(-1.9,0);
     		\draw [dashed, color= black] (1.52,2.25)--(1.52,0);
     		
     		\draw [dashed, color=black] (-1.9,1.9) -- (-1.4,1.95);
     		\draw [thick] (-1.4,1.95) -- (1.5,2.25);
     		
     		\draw [dashed] (-1.9,1.9) -- (-1.46,2.11);
     		\draw [color=black, thick]  (-1.46,2.11) -- (1.9,3.7);

     		\draw[->] (0,0) -- (0,4.2);
     		
     		\draw (1.56,0) node [below] {$b_1^*$};
     		\draw (-1.27,0) node [below] {$b_2^*$};
     		\draw (.8,2.1) node [above] {$\ell_{E[1]}$};
     		\draw (.5,3.) node [above] {$\tilde{\ell}$};
     		\draw (-2.6,1.6) node [above] {$\Pi(E[1])$};
     		\draw (0,0) node [below] {$o$};
     		
     		\draw  (-1.95,0) node [below ] {\small$\mu(E)$};

     		\begin{scriptsize}
     		\fill (0, 0) circle (1.3pt);
     		
     		\fill (-1.9,0) circle (1.3pt);
     		\fill [color=black] (-1.9, 1.9) circle (1.3pt);
     		
     		\fill [color=black] (1.51,2.25) circle (1.3pt);
     		\fill [color=black] (-1.4,1.95) circle (1.3pt);
     		\fill [color=black] (-1.4,0) circle (1.3pt);
     		\fill [color=black] (1.52,0) circle (1.3pt);
     		\fill [color=black] (1.92,3.7) circle (1.3pt);
     		\fill [color=black] (-1.46,2.11) circle (1.3pt);
     		
     		\end{scriptsize}
     		
     		\end{tikzpicture}
     		
     		\caption{The line $\ell_{E[1]}$}
     		
     		\label{fig.line-ell-E[1]}
     		
     	\end{centering}
     	
     \end{figure}
     
   \noindent To prove the claim we only need to show there is no wall for $E[1]$ above the line $\ell_{E[1]}$ to the right of $\Pi(E)$, see Figure \ref{fig.line-ell-E[1]}. Suppose there is such a wall $\tilde{\ell}$ above $\ell_{E[1]}$ whose intersection point with $\partial U$ has $b$-values $\tilde{b}_2 < \tilde{b}_1$ satisfying $\tilde{b}_2 < b_2^*$ and $b_1^* < \tilde{b}_1$. Let $E_1 \hookrightarrow E[1] \twoheadrightarrow E_2$ be a destabilising sequence which makes $E[1]$ unstable below $\tilde{\ell}$. Taking cohomology from the destabilising sequence gives the long exact sequence of coherent sheaves  
	\begin{equation}\label{seqeunce-f}
	0 \rightarrow \cH^{-1}(E_1) \rightarrow E \xrightarrow{f}  \cH^{-1}(E_2) \rightarrow \cH^{0}(E_1) \rightarrow 0. 
	\end{equation}
	We first show that $\cH^{0}(E_1)$ is of rank zero. Suppose it is of rank $r' \neq 0$, then 
	\begin{align*}
	H^{n-1}\ch_1(E) = \ & H^{n-1}\ch_1(\cH^{-1}(E_1)) +H^{n-1}\ch_1(\cH^{-1}(E_2))  - H^{n-1}\ch_1(\cH^{0}(E_1))\\
	< \ & H^n\tilde{b}_2(\ch_0(E) + r') -H^n\tilde{b}_1r' \\
	< \ & H^nb_2^*\ch_0(E)-H^n(b_1^*-b_2^*)
	\end{align*}
	which is not possible by \eqref{def-ell-[1]}. Thus $\cH^{0}(E_1)$ is of rank zero and the exact sequence \eqref{seqeunce-f} gives 
	\begin{equation}\label{slope-min}
	\mu(\im f ) \leq \mu(\cH^{-1}(E_2)). 
	\end{equation} 
	Since $\tilde{\ell}$ crosses the vertical line $b=b_2^*$ at a point inside $U$, we gain
	\begin{equation*}
	0 < \text{Im}[Z_{b_2^*, w}(E_2)] = -H^{n-1}\ch_1(\cH^{-1}(E_2)) + b_2^*H^n\ch_0(\cH^{-1}(E_2)).
	\end{equation*}
	Combining this with \eqref{slope-min} gives
	\begin{equation*}
	b_2^*- \mu(E) > \mu(\cH^{-1}(E_2))-\mu(E) \geq \mu(\im f) -\mu(E) \geq \mu^{\min}(E) -\mu(E) 
	\end{equation*}
	which is not again possible by \eqref{def-ell-[1]}. 
\end{proof}
Proposition \ref{prop.line of stability-E} and \ref{prop.line of stability-E[1]} can be directly shown using \cite[Theorem 1.3 and 1.4]{sun:tilt-stability-vanishing-theorems}. However, we gave here independent proof. We finish this section by mentioning another way to control the position of walls.
\begin{Lem}\label{lem.minimal}
	Take an object $E \in \Coh^{b_0}(X)$ which is $\nu\_{b_0,w_0}$-semistable for some $w_0 > \frac{b_0^2}{2}$. If 
	\begin{equation*}
	0 < \ch_1(E)H^{n-1} -b_0H^n\ch_0(E)\, = \, \min \left\{  	\ch_1(F)H^{n-1} -b_0H^n\ch_0(F) \neq 0 \ | \  F \in \coh^{b_0}(X) \right\},    
	\end{equation*}
then $E$ is $\nu_{b_0,w}$-stable for all $w > \frac{b_0^2}{2}$. 
\end{Lem} 
\begin{proof}
	Assume otherwise, so when we move along the vertical line $b=b_0$, $E$ gets strictly $\nu_{b_0, w}$-semistable for some $w > \frac{b_0^2}{2}$. 
	
	\noindent There is an exact sequence $E_1 \hookrightarrow E \twoheadrightarrow E_2$ in $\Coh^{b_0}(X)$ such that $\nu\_{b_0, w}(E_i) = \nu\_{b_0, w}(E) < +\infty$, thus Im$[Z_{b_0, w}(E_i)] \neq 0$ for $i=1,2$. Since 
	\begin{equation*}
	\text{Im}[Z_{b_0, w}(E_1)] + \text{Im}[Z_{b_0, w}(E_2)] = \text{Im}[Z_{b_0, w}(E)],
	\end{equation*} 
	we gain  Im$[Z_{b_0, w}(E_i)] < \text{Im}[Z_{b_0, w} (E)]$ which is not possible by our assumption.   
\end{proof}

\section{slope stability of the restricted sheaf}

In this section, we prove our main result Theorem \ref{thm.main thm} which provides sufficient conditions for a slope-stable reflexive sheaf whose restriction to a hypersurface is slope semistable. Firstly, we examine $\mu$-stability of a sheaf supported on a hypersurface via the Chern character of its push-forward.  
\begin{Lem}\label{lem.slope-stability}
	Take an irreducible divisor $D \in |k H |$ with the embedding $i \colon D \hookrightarrow X$. A coherent sheaf $F$ on $D$ is  
	$\mu$-(semi)stable if for any quotient sheaf $F \twoheadrightarrow  F'$, we have 
	\begin{equation}\label{slope}
	\dfrac{\ch_2(i_*F).H^{n-2}}{\ch_1(i_*F).H^{n-1}} <(\leq)\ \dfrac{\ch_2(i_*F').H^{n-2}}{\ch_1(i_*F').H^{n-1}}\ .
	\end{equation} 
\end{Lem}
\begin{proof}
	For any coherent sheaf $Q$ on $D$, by adjuction, we get 
	\begin{equation*}
	\chi(\cO_D(-m), Q) = \chi(\cO_X(-m) , i_*Q). 
	\end{equation*}
	Since $X$ is a smooth projective variety, Hirzebruch-Riemann-Roch Formula gives
		\begin{align*}
	\chi(\cO_X(-m) , i_*Q) = \ & \frac{\ch_1(i_*Q)H^{n-1}}{(n-1)!} \ m^{n-1}  \\
	&\ + \left(\frac{\ch_2(i_*Q)H^{n-2}}{(n-2)!} - \frac{\ch_1(i_*Q).K_X.H^{n-2}}{2(n-2)!} \right) m^{n-2} + \dots 
	\end{align*} 
	Thus 
	\begin{equation*}
	\rk(Q) = \frac{\alpha_{n-1}(i_*Q)}{\alpha_{n-1}(i_*\cO_D)} = \frac{\ch_1(i_*Q)H^{n-1}}{k H^n}\ ,  \qquad \deg(Q) = \alpha_{n-2}(i_*Q) - \rk(Q) \ \alpha_{n-2}(i_*\cO_D). 
	\end{equation*}
	Since $Q$ is supported on an irreducible divisor $D \in |k H|$, $\ch_1(i_*Q)$ is a multiple of $H$, so  
	\begin{equation*}
	\deg(Q) = \frac{\ch_2(i_*Q)H^{n-2}}{(n-2)!} + \frac{k\ch_1(i_*Q)H^{n-1}}{2(n-2)!}\,.
	\end{equation*}
	Therefore, if $\alpha_{n-1}(Q) \neq 0$, the slope of $Q$ is 
	\begin{equation*}
	\mu(Q) = \frac{\deg(Q)}{\rk(Q)} = \frac{k H^n}{(n-2)!}\ \frac{\ch_2(i_*Q)H^{n-2}}{\ch_1(i_*Q)H^{n-1}} + \frac{k^2H^n}{2(n-2)!}\,. 
	\end{equation*}
	Hence, for a quotient sheaf $F \twoheadrightarrow F'$ on $D$, we have $\mu(F) < (\leq)\ \mu(F')$ if and only if \eqref{slope} holds.    
\end{proof}
 

The next proposition describes how we can use the 2-dimensional slice $U$ of weak stability conditions 
to prove the $\mu$-stability of the restricted sheaf.    

\begin{Prop}\label{prop.traiangle-stability}
	Let $E$ be a reflexive sheaf on $X$ such that for a point $(b,w) \in U$, $E$ and $E(-kH)[1]$ are $\nu\_{b,w}$-semistable objects in $\coh^b(X)$ with 
	\begin{equation}\label{slope assumption}
	\nu\_{b,w}(E) \leq \nu\_{b,w}(E(-k H)[1]).
	\end{equation}
	Then the restricted sheaf $E|_{D}$ for any irreducible divisor $D \in |k H|$ is $\mu$-semistable on $D$ if it has no quotient sheaf $E|_D \twoheadrightarrow F$ whose point $Z_{b,w}(F)$ lies inside the triangle $\triangle opp'$ or on its sides except $\overline{op'}$ where $p = Z_{b,w}(E)$ and $p' = Z_{b,w}(E|_D)$.  
\end{Prop}
\begin{figure} [h]
	\begin{centering}
		
		\begin{tikzpicture}[line cap=round,line join=round,>=triangle 45,x=1.0cm,y=1.0cm]
		
		
		
		\draw (5,0) node [right] {\small$\text{Re}[Z_{b,w}] =  -\ch_2H^{n-2} +wH^n\ch_0$};
		\draw (3.5,3.2) node [above] {\small $\text{Im}[Z_{b,w}] = \ch_1H^{n-1} -bH^n\ch_0$};
		
		\draw[->,color=black] (-2,0) -- (5,0);
		\draw[->,color=black] (2,0) -- (2,3.2);
		

		\draw [ color=black] (2,0)--(4,1);
		\draw [ color=black] (4,1)--(.1,3);
		\draw [ color=black] (2, 0)--(.1,3);
		
			\draw [ color=black, dashed] (1.8,1.05)--(2, 0);
		\draw [ color=black, dashed] (1.8,1.05)--(.1,3);

		\draw (2.6, .65) node [above] {\small $Z_{b,w}(F)$};
		\draw (4, 1) node [right] {\small$p=Z_{b,w}(E)$};
		\draw (.25, 3.05) node [left] {\small $Z_{b,w}(E|_D)=p'$};
		\draw (2,0) node [below] {$o$};

		\begin{scriptsize}
		\fill (1.8,1.05) circle (1.4pt);
		\fill (2, 0) circle (1.4pt);
		
		\fill (4,1) circle (1.4pt);
		\fill [color=black] (.1, 3) circle (1.4pt);

		\end{scriptsize}
		
		\end{tikzpicture}
		
		\caption{The triangle $\triangle opp'$}
		
		\label{fig.trinagle}
		
	\end{centering}
	
\end{figure}
\begin{proof}
	Since $E$ and $E(-kH)[1]$ are in the heart $\coh^b(X)$, the exact triangle 
	\begin{equation}\label{scc.1}
	E \hookrightarrow E|_D \twoheadrightarrow E(-k H)[1]
	\end{equation}
	implies that $E|_D \in \coh^b(X)$. 
	
	Let $F^{+} \hookrightarrow E|_D$ be the subobject of $E|_D$ of maximum $\nu_{b,w}$-slope in the HN filtration of $E|_D$. If $\nu_{b,w}(F^{+}) > \nu_{b,w}(E(-k)[1])$, then \eqref{slope assumption} implies
	\begin{equation*}
	\text{hom}(F^{+}, E(-k H)[1]) = \text{hom}(F^+, E) = 0
	\end{equation*}
   which is not possible by the exact sequence \eqref{scc.1}, thus the maximum slope $\nu^{\max}_{b,w}(E|_D)$ in the HN filtration of $E|_D$ satisfies 
   \begin{equation}\label{a.1}
   \nu^{\max}_{b,w}(E|_D) \leq \nu_{b,w}(E(-k H)[1]). 
   \end{equation}  
	A similar argument also implies 
	\begin{equation}\label{a.2}
	\nu_{b,w}(E) \leq \nu^{\min}_{b,w}(E|_D)
	\end{equation}
  where $\nu^{\min}_{b,w}(E|_D)$ is the minimum $\nu_{b,w}$-slope in the HN filtration of $E|_D$.  
	
	Suppose $E|_D$ is not $\mu$-semistable on $D$. Then there is a sequence of coherent sheaves 
	\begin{equation}\label{sec}
	F' \hookrightarrow E|_D \twoheadrightarrow F 
	\end{equation}
	such that $\mu(F) < \mu(E|_D) < \mu(F')$. By definition of the heart $\coh^b(X)$, any torsion sheaf $Q$ is in the heart and its $\nu_{b,w}$-slope is equal to $\frac{\ch_2(Q).H^{n-2}}{\ch_1(Q)H^{n-1}}$ if $\ch_1(Q)H^{n-1} \neq 0$, otherwise it is $+\infty$. Thus \eqref{sec} is an exact sequence in $\coh^b(X)$ and by Lemma \ref{lem.slope-stability},
	\begin{equation*}
	\nu\_{b,w}(F)  < \nu\_{b,w}(E|_D) <  \nu\_{b,w}(F').   
	\end{equation*} 
	This implies $Z_{b,w}(F)$ lies to the right of $\overline{op'}$, see Figure \ref{fig.trinagle}. Moreover, we have
	\begin{equation*}
	 \nu_{b,w}^{\min}(E|_D)\leq \nu_{b,w}(F) < \nu_{b,w}(F') \leq  \nu_{b,w}^{\max}(E|_D).   
	\end{equation*}
	Combining this with \eqref{a.1} and \eqref{a.2} shows that
	\begin{equation}\label{slope-order}
	\nu\_{b,w}(E)\leq \nu\_{b,w}(F) \qquad \text{and} \qquad \nu\_{b,w}(F') \leq  \nu\_{b,w}(E(-kH)[1]).    
	\end{equation}
	The first inequality implies $Z_{b,w}(F)$ cannot be to the right of $\overline{op}$, and the second one shows $Z_{b,w}(F)$ cannot be above $\overline{p'p}$. Note that the slope of the line segment $\overline{p'p}$ corresponds to $\nu\_{b,w}(E(-kH)[1])$. Hence $Z_{b,w}(F)$ lies inside $\triangle opp'$ or on the line segments $\overline{op}$ or $\overline{pp'}$. 
\end{proof}

Proposition \ref{prop.traiangle-stability} in particular, implies the following result. 

\begin{Cor}\label{cor-slope-stability}
	Let $E$ be a slope-stable reflexive sheaf such that $E$ and $E(-k H) [1]$ are $\nu\_{b,w}$-stable of the same slope, then $E|_{D}$ is $\mu$-stable for any irreducible divisor $D \in |kH|$. 
\end{Cor}
\begin{proof}
	By proposition \ref{prop.traiangle-stability} $E|_D$ is $\mu$-semistable. Suppose it is strictly $\mu$-semistable and $E|_D \twoheadrightarrow F$ is a proper quotient sheaf with $\mu(E|_D) = \mu(F)$. We know $F$ is also a quotient object of $E|_D$ in $\coh^{b}(X)$ with the same $\nu\_{b, w}$-slope as $E|_D$. Thus $E|_D$ is strictly $\nu\_{b, w}$-semistable. 
	
	Since $E$ is $\nu\_{b, w}$-stable, the composition $E \hookrightarrow E|_D \twoheadrightarrow F$ in $\coh^{b}(X)$ must either be zero or injective. It cannot be zero because this would give a surjection $E(-k H)[1] \twoheadrightarrow F$ in $\coh^{b}(X)$, contradicting the $\nu\_{b, w}$-stability of $E(-k H)[1]$. Thus it is injective, then its cokernel $C$ in $\coh^{b}(X)$ lies in a commutative diagram
	\begin{equation*}
	\xymatrix{
		\,E\,\ar@{^{(}->}[r]\ar@{=}[d] & E|_D \ar@{->>}[r]\ar@{->>}[d]&E(-k H)[1]\ar@{->>}[d]\\
		\,E\,\ar@{^{(}->}[r] &F\ar@{->>}[r]&\,C.}
	\end{equation*}
	Since $E$ and $F$ are $\nu\_{b, w}$-semistable of the same slope, $C$ is also $\nu\_{b, w}$-semistable. Hence the right hand surjection again contradicts the $\nu\_{b,w}$-stability of $E(-k H)[1]$.   
	
\end{proof}

\begin{proof}[Proof of Theorem \ref{thm.main thm}]    	

	Let $\ell$ be the line passing through $\Pi(E)$ and $\Pi(E(-kH))$ for $k >0$, and let $r = H^n\ch_0(E)$. If 
	\begin{equation}\label{c.1}
	 \frac{4}{r^2} \Delta_H(E) < k^2\ ,
	\end{equation}
	then $\ell$ intersects the boundary $\partial U$ at two points with $b$-values $b_2 <b_1$ where    
	\begin{equation*}
	b_1, b_2 = \mu(E) - \dfrac{k}{2} \pm \sqrt{\dfrac{k^2}{4} - \dfrac{\Delta_H(E)}{r^2}}\ .
	\end{equation*}  
	    By Proposition \ref{prop.line of stability-E} and \ref{prop.line of stability-E[1]}, $E$ and $E(-k H)[1]$ are $\nu_{b,w}$-stable for any $(b,w) \in \ell \cap U$ if the following three conditions are satisfied: 
	\begin{enumerate*}
		\item $\frac{r}{H^n}(b_2 - \mu(E(-k H)))= \frac{r}{H^n}(\mu(E) -b_1)< b_1-b_2 = 2\sqrt{\frac{k^2}{4} - \frac{\Delta_H(E)}{r^2}}$
		\item $	\mu(E) -b_1 < \delta(E)$ 
		\item $b_2 -\mu(E(-k H)) < \delta(E)$ 
	\end{enumerate*} 
	where $\delta(E) = \min \left\{ \mu^{\min}(E(-kH))-\mu(E(-kH)) ,\ \mu(E) -\mu^{\max}(E) \right\}$. 
	\begin{figure} [h]
		\begin{centering}
			
			\begin{tikzpicture}[line cap=round,line join=round,>=triangle 45,x=1.0cm,y=1.0cm]
			
			
			\fill [fill=gray!30!white] (0,0) parabola (2,4) parabola [bend at end] (-2,4) parabola [bend at end] (0,0);
			
			\draw (3,0) node [right] {$b, \frac{\ch_1H^{n-1}}{\ch_0H^n}$};
			\draw (0,4.2) node [above] {$w, \frac{\ch_2H^{n-1}}{\ch_0H^n}$};
			
			\draw[->,color=black] (-3,0) -- (3,0);

			\draw [] (0,0) parabola (2,4); 
			\draw [] (0,0) parabola (-2,4);

			\draw [color=black] (-1.6,.6)--(2,2.5); 
			\draw [color=black, dashed ] (-1.6,4)--(-1.6,0);
			\draw [color=black, dashed ] (-.95,.93)--(-.95,0);
			\draw [color=black, dashed] (2,0)--(2,4);
			
			\draw [dashed, color=black] (1.5,2.25)  -- (1.5, 0);

			\draw[->] (0,0) -- (0,4.2);
			
			
			\draw (1.48,0) node [below] {\small $b_1$};
			\draw (-.8,0) node [below] {\small $b_2$};
			\draw (-1.6,.6) node [left] {\small $\Pi(E(-kH))$};
			\draw (2,2.5) node [right] {\small $\Pi(E)$};
			\draw (-1.88,0) node [below] {\small$\mu(E)-k$};
			\draw (2.17,0) node [below] {\small $\mu(E)$};
			\draw (-.6, 1.7) node [below] {$\ell$};
			\begin{scriptsize}
			\fill [color=black] (-0.95,0) circle (1.4pt);
			
			\fill [color=black] (2,0) circle (1.4pt);
			\fill [color=black] (-1.6,0) circle (1.4pt);
			\fill [color=black] (-1.6,.6) circle (1.4pt);
			\fill [color=black] (2,2.5) circle (1.4pt);
			\fill [color=black] (0,0) circle (1.4pt);
			\fill [color=black] (1.5,0) circle (1.4pt);
			\end{scriptsize}
			
			\end{tikzpicture}
			
			\caption{Stability of the restricted bundle}
			
			\label{the line L}
			
		\end{centering}
		
	\end{figure} 

    \noindent The condition (a) is equivalent to  
	\begin{align}\label{c.2}
	\frac{r}{H^n}\left(  \dfrac{k}{2} - \sqrt{\dfrac{k^2}{4} - \dfrac{\Delta_H(E)}{r^2}} \right) < 2\,\sqrt{\dfrac{k^2}{4} - \dfrac{\Delta_H(E)}{r^2}}\,.
	\end{align}
	and thus
	\begin{align}\label{a-1}
	\left(\frac{2H^n}{r}+1\right)^2 \Delta_H(E) < H^n\left(r + H^n\right)k^2.  
	\end{align}
	Note that \eqref{c.2} and so condition (a) implies \eqref{c.1}. Also (b) and (c) are equivalent to  
	\begin{align}\label{eq}
	 \frac{k}{2} - \sqrt{\frac{k^2}{4} - \frac{\Delta_H(E)}{r^2}} < \delta(E).  
	\end{align} 
	Hence, if we have strict inequalities in \eqref{bound for ell-new}, then $E$ and $E(-kH)[1]$ are $\nu\_{b,w}$-stable of the same slope for $(b,w) \in \ell \cap U$, so $E|_D$ is $\mu$-stable by Corollary \ref{cor-slope-stability}.   
	
	Now suppose we have equality in (a), (b), or (c) which means equality in one of the inequalities in \eqref{bound for ell-new}. The computation in \eqref{c.2} shows that still \eqref{c.1} holds, so $\ell \cap U \neq \emptyset$. Then the structure of walls and Proposition \ref{prop.line of stability-E} and \ref{prop.line of stability-E[1]} imply that $E$ and $E(-kH)[1]$ are $\nu_{b,w}$-semistable for $(b,w) \in \ell \cap U$. Hence $E|_D$ is $\mu$-semistable on $D$ by Proposition \ref{prop.traiangle-stability}.      	
\end{proof}
As a consequence of Theorem \ref{thm.main thm}, we obtain a variant of Langer's restriction Theorem \cite[Theorem 5.2]{langer:semistable-sheaves-in-positive-characteristics}.   
\begin{Cor}\label{cor-langer}
	Let $E$ be a $\mu$-stable reflexive sheaf of rank $\rk >1$. The restricted sheaf $E|_D$ for any irreducible divisor $D \in |kH|$ is $\mu$-stable if 
\begin{equation}\label{l}
k>H^n\rk\big(\rk -1\big)\widetilde{\Delta}(E)+ \frac{1}{H^n\rk\big(\rk -1\big)} \,.  
\end{equation}
\end{Cor}
\begin{proof}
We always have 
	\begin{equation*}
	\delta(E) \geq \frac{1}{H^n\rk(\rk-1)}\,.
	\end{equation*}
  Hence if \eqref{l} holds, then the second inequality in \eqref{bound for ell-new} is satisfied.        
Also, if $\rk \geq 3$ or $\Delta_H(E) \geq 3$, then  
 \begin{align}\label{co}
 \frac{\left(\rk+2\right)^2}{\rk +1} \frac{\Delta_H(E)}{(H^n\rk)^2} \leq \left( \frac{1}{\rk(\rk-1)} + H^n\rk(\rk-1)\frac{\Delta_H(E)}{ (H^n\rk)^2}\right)^2 < k^2, 
 \end{align}
which shows the first inequality in \eqref{bound for ell-new} holds, so the claim follows.

If $\rk =2$ and $0 \leq \Delta_H(E) \leq 2$, then the assumption \eqref{l} implies that \eqref{c.1} holds. Thus there are points $(b_0,w_0) \in U$ where $E$ and $E(-kH)[1]$ are in $\coh^{b_0}(X)$ and have the same $\nu_{b_0,w_0}$-slope. Moreover, the same argument as in Step 2 in the proof of Theorem \ref{thm.bounding giesker chamber} implies that there are no wall for $E$ and $E(-kH)[1]$ in $U$. Thus $E$ and $E(-kH)[1]$ are $\nu_{b_0, w_0}$-stable and so the claim follows from Corollary \ref{cor-slope-stability}.        
\end{proof}
\begin{Rem}\label{rem-langer}
\cite[Theorem 5.2]{langer:semistable-sheaves-in-positive-characteristics} covers more general cases than Corollary \ref{cor-langer}, it allows $X$ to be a smooth projective variety over an arbitrary algebraically closed field, and assumes $E$ to be a $\mu$-stable torsion free sheaf with respect to some nef divisors.
\end{Rem}

In some special cases, we have more control over the position of the walls for $E$ and $E(-kH)[1]$, so we gain stronger effective restriction results. 
For instance, suppose $(X,H)$ is a smooth polarised complex projective variety of dimension $n \geq 2$ such that 
\begin{equation}\label{ass}
\text{$H^n$ divides $H'.H^{n-1}$ for all $H' \in \Pic(X)$.}
\end{equation}
Let $E$ be a $\mu$-stable reflexive sheaf with $\gcd\left(H^n\ch_0(E) , \ch_1(E)H^{n-1}\right) =H^n$. We define 
	\begin{equation*}
m(E)^{\pm} = \min \left\{m \in \mathbb{N}\colon \ \  m' \,H^n\ch_0(E)\, - m\,  \ch_1(E)H^{n-1}\, = \pm H^n  \ \ \text{for some} \ m' \in \mathbb{Z} \right\}. 
\end{equation*}  

\begin{Prop}\label{special}
	Let $E$ be a $\mu$-stable reflexive sheaf as above with rank $\rk > 1$. 
	The restricted sheaf $E|_D$ for any irreducible divisor $D \in |kH|$ is $\mu$-(semi)stable on $D$ if 
	\begin{equation}\label{a}
	k >\! (\geq)\ \,  \frac{\widetilde{\Delta}(E)}{\delta} + \delta\,,
	\end{equation}
	where $\delta = \frac{1}{\rk} \times\min\{\frac{1}{m(E)^+} , \frac{1}{m(E)^-}\}$. In particular, $E|_D$ is $\mu$-(semi)stable if
	\begin{equation}\label{f}
k > \!(\geq) \, \rk(\rk -1) \widetilde{\Delta}(E) + \frac{1}{\rk(\rk -1)} \, . 
\end{equation} 
\end{Prop} 
\begin{proof}
	First of all, since $\delta \leq \frac{1}{2}$, \eqref{a} is equivalent to 
	\begin{equation}\label{a.11}
	\frac{k}{2} - \sqrt{\frac{k^2}{4} - \widetilde{\Delta}(E)} \ <\!(\leq)\ \delta\,.
	\end{equation}	
	As in the proof of Theorem \ref{thm.main thm}, let $\ell$ be the line passing through $\Pi(E)$ and $\Pi(E(-kH))$. First assume the inequality in \eqref{a} is strict, so 
	\begin{equation}\label{1}
	k^2 - 4 \widetilde{\Delta}(E)  > \left(\frac{\widetilde{\Delta}(E)}{\delta} - \delta \right)^2. 
	\end{equation} 
	Thus \eqref{c.1} holds and so $\ell \cap U \neq \emptyset$. Let $b_2 <b_1$ be the values of $b$ at the intersections points of $\ell$ with $\partial U$. Hence \eqref{a.11} with applying the same argument as in the proof of Theorem \ref{thm.main thm} implies
	\begin{equation}\label{2}
	\mu(E) -b_1 < \delta \qquad \text{and} \qquad b_2 -\mu(E(-kH)) < \delta. 
	\end{equation}
	 Therefore the line $\ell$ intersects the vertical lines $b = \mu(E) -\frac{1}{\rk m(E)^-}$ and $b = \mu(E(-kH)) +\frac{1}{\rk m(E)^+}$ at points inside $U$. So by Lemma \ref{lem.minimal}, $E$ and $E(-kH)[1]$ are $\nu_{b,w}$-stable for any $(b,w) \in \ell \cap U$ and the claim follows from Corollary \ref{cor-slope-stability}. If we have equality in \eqref{a} which means equality in \eqref{1} or \eqref{2}, we obtain the claim by applying Proposition \ref{prop.traiangle-stability} at the limiting points $(b,w) \in U$ just above $\ell$.    
	 
    Finally, by definition of $m(E)^{\pm}$, we know $\delta \geq \frac{1}{\rk(\rk -1)}$. Thus \eqref{a.11} implies $E|_D$ for $D \in |kH|$ is $\mu$-(semi)stable if   
	 \begin{equation*}
	 \frac{k}{2} - \sqrt{\frac{k^2}{4} - \widetilde{\Delta}(E)} <\! (\leq)\ \frac{1}{\rk(\rk -1)}
	 \end{equation*}	
	 which is equivalent to \eqref{f}.  
\end{proof}

\section{Clifford indices of curves over K3 surfaces }\label{section.clifford}
In this section, we assume $(S,H)$ is a smooth polarized K3 surface over $\mathbb{C}$ such that 
\begin{equation*}
\text{$H^2$ divides $H.D$ for all curve classes $D$ on $S$}. 
\end{equation*}
Let $\iota \colon C \hookrightarrow S$ be any smooth curve in the linear system $|H|$ of genus $g$. It is proved in \cite{lazarsfeld:brill-noether} and \cite[Theorem 1.1]{bayer:brill-noether} that there exists a globally generated degree $d$ line bundle $A$ on $C$ with $r$ global sections if and only if 
\begin{equation}\label{rho-1}
\rho(r-1,d,g) = g- r(g-1-d+r) \geq 0\ . 
\end{equation} 
Take such a line bundle $A$ on $C$. Then the kernel of the evaluation of sections of $\iota_*A$
\begin{equation}\label{sec.1}
0 \rightarrow F_{C,A} \rightarrow H^0(C,A) \otimes \cO_S \xrightarrow{\text{ev}} \iota_*A \rightarrow 0
\end{equation}
is a vector bundle of rank $r$. The Lazarsfeld-Mukai bundle $E_{C,A}$ associated to the pair $(C, A)$ is the dual of $F_{C,A}$ and it is of character 
\begin{equation*}
\ch(E_{C,A}) =  \big(\, r,\, H,\, g-1-d\, \big). 
\end{equation*}  
The bundle $E_{C,A}$ has been appeared, for example, in Lazarsfeld's proof of Brill-Noether Petri Theorem \cite{lazarsfeld:brill-noether}, in the Mukai's classification of prime Fano manifolds of coindex 3, or in Voisin’s proof of Green’s canonical syzygy conjecture \cite{voisin:greens-conjecture}; see \cite{aprodu:lazarsfeld} for a survey of applications.

The $\mu$-slope of any destabilising subsheaf of $F_{C,A}$ must be less than zero because of the exact sequence \eqref{sec.1}, thus $F_{C,A}$ and so $E_{C,A}$ are $\mu$-stable. Hence \eqref{rho-1} implies that for any $\ch_2 \in \mathbb{Z}$ with 
\begin{equation}\label{bundle}
0 \leq \ \ch_2 \ \leq \frac{g}{r} -r\,,
\end{equation}
there is a stable vector bundle $E_{C,A}$ of Chern character $v = (r, H, \ch_2)$ for any curve $C \in |H|$. 

The restriction of Lazarsfeld-Mukai bundle $E_{C,A}$ to curves on the K3 surface $S$ has led to counterexamples to ${\bf M_3}(C)$ \cite{farkas:higher-rank-brill-noether} and ${\bf M_4}(C)$ \cite{aprodu:restricted-lazarsfeld-mukai}. By applying Theorem \ref{thm.main thm}, we extend these results to higher ranks.  

\begin{proof}[Proof of Theorem \ref{thm.k3 surface.cor}]
By Proposition \ref{special}, inequality \eqref{f}, we know $E|_C$ is $\mu$-semistable if 
\begin{align*}
1 \geq\ & r(r-1) \widetilde{\Delta}(E) + \frac{1}{r(r-1)} \\
 \geq \ &\frac{r-1}{rH^2} \left(H^2 -2r\ch_2\right) + \frac{1}{r(r-1)}
\end{align*}
which clearly holds when $\ch_2 \geq 0$. 
Since $E$ is $\mu$-stable of positive slope, hom$(E, \cO_X) = 0$, thus 
\begin{equation*}
h^0(S, E) = \chi(\cO_S, E) + h^1(S, E) \geq \chi(\cO_S, E) = 2r+\ch_2. 
\end{equation*} 
Hence if $r \geq 2$ and $\ch_2 \geq 0$, the restricted bundle $E|_C$ is $\mu$-semistable with
\begin{equation*}
\deg(E|_C)= 2g-2 \leq r(g-1) \qquad \text{and} \qquad 2r \leq h^0(S, E) \leq  h^0(C, E|_C)\ .   
\end{equation*}       
The last inequality follows from hom$(\cO_S, E(-H)) = 0$. So $E|_C$ contributes to the rank $r$-Clifford index of $C$. 

By \eqref{rho-1}, we obtain $\Cliff(C) = \left\lfloor \frac{g-1}{2} \right\rfloor$. Moreover, 
\begin{equation*}
\Cliff(E|_C) \leq \frac{2g-2}{r} - \frac{2}{r} (r+\ch_2)   
\end{equation*}   
and the right hand side is clearly less than $\left\lfloor \frac{g-1}{2} \right\rfloor$ when (i) $r \geq 4$, or  (ii) $r=3$ and 
	\begin{equation}\label{less}
	g -\ch_2 < 4+ \frac{3}{2}\left\lfloor \frac{g-1}{2} \right\rfloor. 
	\end{equation}
If $r \geq 4$ and $g \geq r^2$, then by \eqref{bundle}, there is a stable vector bundle $E$ on $S$ with Chern character $(r, H, \ch_2)$ such that $\ch_2 \geq 0$. The above argument implies $E|_C$ contributes to the rank $r$-Clifford index of $C$ and 
\begin{equation*}
\Cliff_r(C) \leq \Cliff(E|_C) <  \left\lfloor \frac{g-1}{2} \right\rfloor, 
\end{equation*}    
thus ${\bf M_r}(C)$ fails. If $r=3$, then by \eqref{bundle} there is a stable vector bundle $E$ of class $(3, H, \ch_2)$ when $\ch_2 =  \lfloor \frac{g}{3} \rfloor -3$. Then \eqref{less} implies $\Cliff(E|_C)$ is less than the expected number if
\begin{equation*}
g - \left\lfloor \frac{g}{3} \right\rfloor + 3 < 4+ \frac{3}{2}\left\lfloor \frac{g-1}{2} \right\rfloor
\end{equation*}
which holds for $g=7, 9$ and $g \geq 11$ as claimed.      
\end{proof}
\section{slope-stability of the restriction of tangent bundle of $\mathbb{P}^n$}\label{sec3}
Let $S$ be a smooth projective K3 surface and let $L$ be an ample line bundle generated by its global sections. The kernel of the evaluation of sections of $L$ is denoted by $M_L$ \eqref{def M_L}. In this section, we prove $M_{L}$ is $\mu$-stable with respect to $L$. 

Let $L =kH$ where $H$ is a primitive ample line bundle on $S$. Similar to section \ref{section.2.background}, we define $\nu_{b,w}$-stability and the projection $\Pi$ \eqref{pi} with respect to the polarisation $H$. On K3 surfaces, we can slightly enlarge the $2$-dimensional slice of stability conditions $U$ due to a stronger Bogomolov-type inequality. Recall that an object $E \in \mathcal{D}(S)$ is called spherical if $\text{Hom}_{\mathcal{D}(S)}(E,E[i]) = \mathbb{C}$ for $i=0,2$ and it is zero otherwise.   
\begin{Lem}\label{lem-spece-extended}
	For any point $(b,w) \in \mathbb{R}^2$ with $w > \frac{b^2}{2} - \frac{1}{H^2}$, the pair $\left(\coh^b(X), Z_{b,w}\right)$ is a Bridgeland stability condition on $S$ if for any spherical bundle $E$ on $S$ with $\mu(E) =b$, we have $\frac{\ch_2(E)}{H^2\ch_0(E)} < w$.   
\end{Lem}
\begin{proof}
	As it is shown in the proof of \cite[Lemma 6.2]{bridgeland:K3-surfaces}, we only need to check that for any slope-stable vector bundle $E$ on $S$ with Im$[Z_{b,w}(E)] = 0$, we have Re$[Z_{b,w}(E)] > 0$. We know Hom$(E, E) = \text{Hom}(E, E[2]) = \mathbb{C}$, thus Riemann-Roch theorem gives
	\begin{align*}
	-\chi(E,E) = \ & \sum_{i}(-1)^{i+1} \dim _{\mathbb{C}} \text{Hom}_X^i(E,E) = -2+ \dim _{\mathbb{C}} \text{Hom}_X(E,E[1])   \\
	 = \ & \big(\ch_1(E)\big)^2-2\ch_0(E)\big(\ch_2(E) + \ch_0(E) \big). 
	\end{align*}    
which implies $-\chi(E, E) \geq -2$. If $E$ is not spherical, i.e. $-\chi(E, E) \geq 0$, then Hodeg index theorem gives 
\begin{equation*}
0 \leq \frac{1}{2}\left(\frac{\ch_1(E).H}{H^2\ch_0(E)} \right)^2 - \frac{\ch_2(E)}{H^2\ch_0(E)} - \frac{1}{H^2} . 
\end{equation*}
Hence if $\mu(E) = b$, we get 
\begin{equation*}
\frac{1}{H^2\ch_0(E)}  \text{Re}[Z_{b,w}(E)] = -\frac{\ch_2(E)}{H^2\ch_0(E)} + w >0,
\end{equation*}
because $w > \frac{b^2}{2} - \frac{1}{H^2}$ by our assumption. Thus we only need to check positivity of the real part for spherical bundles as mentioned.  
\end{proof}
The next step is to control the position of the projection of spherical bundles. Since $h^0(L) = \chi(\cO_X, L) = \frac{L^2}{2} +2$, we get 
\begin{equation*}
\ch(M_{L}) = \left(\frac{L^2}{2} +1 , \ -L, \ -\frac{L^2}{2} \right). 
\end{equation*}
Define 
\begin{equation}\label{v}
V \coloneqq \left\{(b,w) \in \mathbb{R}^2 \colon -\frac{2L.H}{H^2(L^2+2)} \,\leq b <\frac{L.H}{H^2}\ , \ \ \frac{L^2}{2LH}b < w \   \right\}. 
\end{equation}

\begin{figure} [h]
	\begin{centering}
		
		\begin{tikzpicture}[line cap=round,line join=round,>=triangle 45,x=1.0cm,y=1.0cm]
		
		\filldraw[fill=gray!40!white, draw=white] (1.95,5.3)--(1.95,3.8)--(-.23,-.45)--(-.23,5.3)--(1.95, 5.3);
		
		
		\draw[->,color=black] (-3,0) -- (3,0);
			\draw[->] (0,-1) -- (0,5.5);
		
		\draw [] (0,0) parabola (2.3,5.3); 
		\draw [] (0,0) parabola (-2.3,5.3); 
		\draw [color=black, dashed] (-.23,-.45)--(1.95,3.8); 
		
		\draw [color=black] (-.23,-.45)--(-0.23,5.3);
		\draw [color=black, dashed ] (1.95,0)--(1.95,5.3);

		

		
		\draw (-.73,-.45) node [below] {$\Pi(M_{L})$};
		\draw (.1, .08) node [below] {$o$};
		\draw (3.05,3.75) node [left] {$\Pi(L)$};
		\draw (2.45,5.2) node [above] {$w= \frac{b^2}{2}$};
		
		\draw (0,5.5) node [above] {$w, \frac{\ch_2.H}{\ch_0H^3}$};
		\draw (3,0) node [right] {$b, \frac{\ch_1.H^2}{\ch_0H^3}$};
		\draw (1,4) node [above] {$V$};
		\draw (1.95,0) node [below] {$l$};
		\begin{scriptsize}
		
		\fill [color=black] (-.23,-.45) circle (1.4pt);
		\fill [color=black] (1.95,3.8) circle (1.4pt);
		\fill [color=black] (1.95, 0) circle (1.4pt);
		\fill [color=black] (0,0) circle (1.4pt);
		\end{scriptsize}
		
		\end{tikzpicture}
		
		\caption{no projection of spherical objects in the grey area $V$}
		
		\label{the area v}
		
	\end{centering}
	
\end{figure}   
\begin{Lem}\label{lem. no soherial in V}
	If $k = L.H/H^2 > 1$, there is no spherical bundle $E$ on $S$ such that its projection $\Pi(E)$ lies in the area $V$, see Figure \ref{the area v}. 
\end{Lem}  
\begin{proof}
	Suppose for a contradiction that there is such an spherical bundle $E$. Let 
	\begin{equation*}
	(r, c, s) \coloneqq \big(\ch_0(E)H^2,\ \ch_1(E)H,\ \ch_2(E)\big).
	\end{equation*} 
	 Since $\chi(E, E) = 2$, the Hodge index theorem gives
	\begin{align}\label{chi}
	-\frac{1}{H^2(\ch_0)^2} =\  &\frac{(\ch_1)^2}{2H^2(\ch_0)^2} - \frac{\ch_2}{H^2\ch_0} - \frac{1}{H^2}\nonumber\\
	\leq\  & \frac{1}{2}\left(\frac{H\ch_1}{H^2\ch_0} \right)^2 - \frac{\ch_2}{H^2\ch_0} - \frac{1}{H^2} 
	\end{align}
	Thus $\Pi(E)$ cannot be in the area $U$ above $w = \frac{b^2}{2}$, so 
	\begin{equation}\label{position}
	-\frac{2L.H}{H^2(L^2+2)}\leq \frac{H\ch_1}{H^2\ch_0} < 0 \qquad \text{and} \qquad-\frac{L^2}{H^2(L^2+2)} < \frac{\ch_2}{H^2\ch_0}
	\end{equation}
  Combining \eqref{chi} and \eqref{position} gives 
  \begin{align*}
  -\frac{L^2}{H^2(L^2+2)} 
  < \frac{\ch_2}{H^2\ch_0} & \ \leq  \frac{1}{2}\left(\frac{H\ch_1}{H^2\ch_0} \right)^2 - \frac{1}{H^2} +	\frac{1}{H^2(\ch_0)^2}\\
  & \ \leq \frac{1}{2}\left(\frac{2L.H}{H^2(L^2+2)}\right)^2 - \frac{1}{H^2} +	\frac{1}{H^2(\ch_0)^2}
  \end{align*}
  The first inequality gives $-1 < \frac{\ch_2}{\ch_0}$ and comparing the first and the last sentences gives $\ch_0 \leq \frac{L^2}{2}+1$. Hence $-\ch_0 +1 \leq \ch_2$ and so
  \begin{align*}
  \frac{1}{\ch_0} - \frac{1}{(\ch_0)^2} \ \leq\ \frac{\ch_2}{\ch_0} +1 -	\frac{1}{(\ch_0)^2} \ \leq\  \frac{H^2}{2}\left(\frac{H\ch_1}{H^2\ch_0} \right)^2 
  \end{align*}
  Here the second inequality comes from \eqref{chi}.  
  If $\ch_0 > 1$, the left hand side is minimum when $\ch_0$ is the maximum value $\frac{L^2}{2} +1$, so we must have 
  \begin{equation*}
  \frac{2L^2}{(L^2+2)^2} \leq \frac{H^2}{2}\left(\frac{H\ch_1}{H^2\ch_0} \right)^2 
  \end{equation*} 
   which is not possible by \eqref{position}. Thus $\ch_0 =1$, so \eqref{position} gives 
   \begin{equation*}
   - \frac{2}{k} <  - \frac{2L.H}{L^2 +2} < H\ch_1 < 0 
   \end{equation*}
which is not again possible for $k> 1$.   
\end{proof}

Therefore, by Lemma \ref{lem-spece-extended}, the pair $\big(\Coh^b(X), Z_{b,w} \big)$ is a Bridgeland stability conditions on $\mathcal{D}^b(X)$ for any $(b,w) \in V$ if $k >1$. 

\begin{proof}[Proof of Theorem \ref{thm. Theorem of ML}]
The line bundle $L$ is in the heart $\Coh^b(X)$ for $b<k = \frac{L.H}{H^2}$. Since $\Delta_H(L) = 0$, it is $\nu_{b,w}$-stable for any $(b,w) \in U$ \cite[Corollary 3.11(a)]{bayer:the-space-of-stability-conditions-on-abelian-threefolds}. Similarly, $\cO_X[1] \in \Coh^b(X)$ for $b \geq 0$ and it is $\nu_{b,w}$-stable for any $(b,w) \in U$. Let $\ell$ be the line passing through $\Pi(L)$ and $\Pi(\cO_X)$. Then $L$ and $\cO_X[1]$ are $\nu_{b_0, w_0}$-stable of the same slope for any $(b_0,w_0) \in \ell \cap U$. Thus the short exact sequence    
\begin{equation*}
L \hookrightarrow M_{L}[1] \twoheadrightarrow \mathcal{O}_X^{\oplus h^0(L)}[1]
\end{equation*}  
in $\Coh^{b_0}(X)$ implies that $M_L[1]$ is also $\nu_{b_0, w_0}$-semistable. 

We claim $M_L[1]$ is $\nu_{b_0, w^+}$-stable where $w_0 < w^+ \ll w_0 +1$. Otherwise, the structure of locally finite set of walls described in Proposition \ref{prop. locally finite set of walls} shows that there is a destabilising sequence 
\begin{equation*}
E_1 \hookrightarrow M_L[1] \twoheadrightarrow E_2
\end{equation*}
such that the $E_i$ have the same $\nu_{b_0, w_0}$-slope as $M_L[1]$ and they make $M_L[1]$ unstable above $\ell$, i.e.
\begin{equation}\label{order-1}
\nu\_{b_0, w}(E_1) \geq \nu\_{b_0, w}(M_L[1]) \geq \nu\_{b_0, w}(E_2)
\end{equation}
for any $w > w_0$. Hence, rank of $E_1$ is negative. We know $\nu_{b_0, w_0}$-stable factors of $E_1$ and $E_2$ are the Jordan-H$\ddot{\text{o}}$lder factors of $M_L[1]$. Since $L$ and $h^0(L)$-copies of $\cO_X[1]$ are the $\nu_{b_0, w_0}$-Jordan-H$\ddot{\text{o}}$lder factors of $M_L[1]$, we gain the stable factors of $E_1$ are either 
\begin{enumerate}
	\item $L$ and $m$-copies of $\cO_X[1]$, or
	\item $m'$-copies of $\cO_X[1]$.
\end{enumerate}
In the first case, $E_2$ must be the direct sum of copies of $\cO_X[1]$ and clearly \eqref{order-1} does not hold. In case (b), we get hom$(\cO_X, M_L) \neq 0$ which is not possible by the definition of $M_L$ \eqref{def M_L}, thus there is no such destabilising sequence and the claim follows.

If $k = \frac{L.H}{H^2} =1$, Lemma \ref{lem.minimal} implies that there is no wall for $M_L[1]$ crossing the vertical line $b = 0$. Thus $\nu_{b_0, w^+}$-stability of $M_L[1]$ implies that it is $\nu_{b =0, w}$-stable for $w \gg 0$. Hence $M_L$ is $\mu$-stable \cite[Proposition 14.2]{bridgeland:K3-surfaces}. 

If $k>1$, we claim there is no wall for $M_{L}[1]$ in the region $V$ \eqref{v} and it is $\nu_{b,w}$-stable for any $(b,w) \in V$. Suppose there is a wall $\ell'$ in $V$ with a destabilising sequence $E_1 \hookrightarrow M_L[1] \twoheadrightarrow E_2$. We know this wall ends at $\Pi(M_L)$. For any $(b,w) \in \ell' \cap V$, we have $Z_{b,w}(E_i) \neq 0$ and   
\begin{equation*}
\abs{Z_{b,w}(E_1)} + \abs{Z_{b,w}(E_2)}  = \abs{Z_{b,w}(M_L[1])}  
\end{equation*}
which implies $\abs{Z_{b,w}(E_i)} < \abs{Z_{b,w}(M_L)}$. If we move $(b,w)$ along the wall $\ell'$ towards $\Pi(M_L)$, then $ \abs{Z_{b,w}(M_L[1])} \rightarrow 0$, thus 
\begin{equation*}
\abs{Z_{b,w}(E_i)} = \abs{\big(\ch_1(E_i)H - b H^2\ch_0(E_i)\big) i + H^2\ch_0(E_i) w - \ch_2(E_i)} \rightarrow 0. 
\end{equation*}
Therefore $\ch_0(E_i) \neq 0$ and $\Pi(E_i) = \Pi(M_{L})$, so the $E_i$ have the same $\nu_{b,w}$-slope as $M_L[1]$ for any $(b,w) \in V$ and cannot make a wall. Thus, in particular, $M_L[1]$ is $\nu_{b, w}$-stable for any $(b,w) \in V$ with $b= \mu(M_L)$, so it is $\mu$-semistable \cite[Lemma 2.7]{bayer:the-space-of-stability-conditions-on-abelian-threefolds}. If $M_L$ is strictly $\mu$-semistable, there is an exact sequence $F_1 \hookrightarrow M_L \twoheadrightarrow F_2$ of $\mu$-semistable coherent sheaves of the same slope. Thus the $F_i[1]$ have the same $\nu_{b=\mu(M_L), w}$-slope as $M_L[1]$ and they are in the heart $\coh^{\mu(M_L)}(X)$. Thus $M_L[1]$ is strictly $\nu_{b=\mu(M_L), w}$-semistable which is not possible by the above argument.     
	
\end{proof}

\bibliography{mybib}
\bibliographystyle{halpha}

\end{document}